\magnification=1200
\catcode`\@=11 
 
\def\nolabels{\def\wrlabel##1{}\def\eqlabel##1{}\def\reflabel##1{}}
\def\writelabels{\def\wrlabel##1{\leavevmode\vadjust{\rlap{\smash%
{\line{{\escapechar=` \hfill\rlap{\sevenrm\hskip.03in\string##1}}}}}}}%
\def\eqlabel##1{{\escapechar-1\rlap{\sevenrm\hskip.05in\string##1}}}%
\def\thlabel##1{{\escapechar-1\rlap{\sevenrm\hskip.05in\string##1}}}%
\def\reflabel##1{\noexpand\llap{\noexpand\sevenrm\string\string\string##1}}}
\nolabels
\global\newcount\secno \global\secno=0
\global\newcount\meqno \global\meqno=1
\global\newcount\mthno \global\mthno=1
\global\newcount\mexno \global\mexno=1
\global\newcount\mquno \global\mquno=1
\global\newcount\tblno \global\tblno=1
\def\newsec#1{\global\advance\secno by1 
\global\subsecno=0\xdef\secsym{\the\secno.}\global\meqno=1\global\mthno=1
\global\mexno=1\global\mquno=1\global\figno=1\global\tblno=1

\bigbreak\medskip\noindent{\bf\the\secno. #1}\writetoca{{\secsym} {#1}}
\par\nobreak\medskip\nobreak}
\xdef\secsym{}
\global\newcount\subsecno \global\subsecno=0
\def\subsec#1{\global\advance\subsecno by1 \global\subsubsecno=0
\xdef\subsecsym{\the\subsecno.}
\bigbreak\noindent{\bf\secsym\the\subsecno. #1}\writetoca{\string\quad
{\secsym\the\subsecno.} {#1}}\par\nobreak\medskip\nobreak}
\xdef\subsecsym{}
\global\newcount\subsubsecno \global\subsubsecno=0
\def\subsubsec#1{\global\advance\subsubsecno by1
\bigbreak\noindent{\it\secsym\the\subsecno.\the\subsubsecno.
                                   #1}\writetoca{\string\quad
{\the\secno.\the\subsecno.\the\subsubsecno.} {#1}}\par\nobreak\medskip\nobreak}
\global\newcount\appsubsecno \global\appsubsecno=0
\def\appsubsec#1{\global\advance\appsubsecno by1 \global\subsubsecno=0
\xdef\appsubsecsym{\the\appsubsecno.}
\bigbreak\noindent{\it\secsym\the\appsubsecno. #1}\writetoca{\string\quad
{\secsym\the\appsubsecno.} {#1}}\par\nobreak\medskip\nobreak}
\xdef\appsubsecsym{}
\def\appendix#1#2{\global\meqno=1\global\mthno=1\global\mexno=1
\global\figno=1\global\tblno=1
\global\subsecno=0\global\subsubsecno=0
\global\appsubsecno=0
\xdef\appname{#1}
\xdef\secsym{\hbox{#1.}}
\bigbreak\bigskip\noindent{\bf Appendix #1. #2}
\writetoca{Appendix {#1.} {#2}}\par\nobreak\medskip\nobreak}
%
%
\def\eqnn#1{\xdef #1{(\secsym\the\meqno)}\writedef{#1\leftbracket#1}%
\global\advance\meqno by1\wrlabel#1}
\def\eqna#1{\xdef #1##1{\hbox{$(\secsym\the\meqno##1)$}}
\writedef{#1\numbersign1\leftbracket#1{\numbersign1}}%
\global\advance\meqno by1\wrlabel{#1$\{\}$}}
\def\eqn#1#2{\xdef #1{(\secsym\the\meqno)}\writedef{#1\leftbracket#1}%
\global\advance\meqno by1$$#2\eqno#1\eqlabel#1$$}
%
%
\def\thm#1{\xdef #1{\secsym\the\mthno}\writedef{#1\leftbracket#1}%
\global\advance\mthno by1\wrlabel#1}
\def\exm#1{\xdef #1{\secsym\the\mexno}\writedef{#1\leftbracket#1}%
\global\advance\mexno by1\wrlabel#1}
%
%
\def\tbl#1{\xdef #1{\secsym\the\tblno}\writedef{#1\leftbracket#1}%
\global\advance\tblno by1\wrlabel#1}
%
\newskip\footskip\footskip14pt plus 1pt minus 1pt 
\def\f@@t{\baselineskip\footskip\bgroup\aftergroup\@foot\let\next}
\setbox\strutbox=\hbox{\vrule height9.5pt depth4.5pt width0pt}
\global\newcount\ftno \global\ftno=0
\def\foot{\global\advance\ftno by1\footnote{$^{\the\ftno}$}}
%
\newwrite\ftfile
\def\footend{\def\foot{\global\advance\ftno by1\chardef\wfile=\ftfile
$^{\the\ftno}$\ifnum\ftno=1\immediate\openout\ftfile=foots.tmp\fi%
\immediate\write\ftfile{\noexpand\smallskip%
\noexpand\item{f\the\ftno:\ }\pctsign}\findarg}%
\def\footatend{\vfill\eject\immediate\closeout\ftfile{\parindent=20pt
\centerline{\bf Footnotes}\nobreak\bigskip\input foots.tmp }}}
\def\footatend{}
%
%
\global\newcount\refno \global\refno=1
\newwrite\rfile
\def\ref{\the\refno\nref}
\def\bref{\nref}
\def\nref#1{\xdef#1{\the\refno}\writedef{#1\leftbracket#1}%
\ifnum\refno=1\immediate\openout\rfile=refs.tmp\fi
\global\advance\refno by1\chardef\wfile=\rfile\immediate
\write\rfile{\noexpand\item{[#1]\ }\reflabel{#1\hskip.31in}\pctsign}\findarg}
\def\findarg#1#{\begingroup\obeylines\newlinechar=`\^^M\pass@rg}
{\obeylines\gdef\pass@rg#1{\writ@line\relax #1^^M\hbox{}^^M}%
\gdef\writ@line#1^^M{\expandafter\toks0\expandafter{\striprel@x #1}%
\edef\next{\the\toks0}\ifx\next\em@rk\let\next=\endgroup\else\ifx\next\empty%
\else\immediate\write\wfile{\the\toks0}\fi\let\next=\writ@line\fi\next\relax}}
\def\striprel@x#1{} \def\em@rk{\hbox{}}
\def\lref{\begingroup\obeylines\lr@f}
\def\lr@f#1#2{\gdef#1{\ref#1{#2}}\endgroup\unskip}

\def\addref#1{\immediate\write\rfile{\noexpand\item{}#1}} 
\def\startrefs#1{\immediate\openout\rfile=refs.tmp\refno=#1}
\def\xref{\expandafter\xr@f}\def\xr@f[#1]{#1}
\def\refs#1{[\r@fs #1{\hbox{}}]}
\def\r@fs#1{\edef\next{#1}\ifx\next\em@rk\def\next{}\else
\ifx\next#1\xref #1\else#1\fi\let\next=\r@fs\fi\next}
%

%
 \newwrite\ffile\global\newcount\figno \global\figno=1
%
%
\def\fig{\the\figno\nfig}
\def\nfig#1{\xdef#1{\secsym\the\figno}%
\writedef{#1\leftbracket \noexpand~\the\figno}%
\ifnum\figno=1\immediate\openout\ffile=figs.tmp\fi\chardef\wfile=\ffile%
\immediate\write\ffile{\noexpand\medskip\noexpand\item{Figure\ \the\figno. }
\reflabel{#1\hskip.55in}\pctsign}\global\advance\figno by1\findarg}
\def\xfig{\expandafter\xf@g}\def\xf@g \penalty\@M\ {}
\def\figs#1{figs.~\f@gs #1{\hbox{}}}
\def\f@gs#1{\edef\next{#1}\ifx\next\em@rk\def\next{}\else
\ifx\next#1\xfig #1\else#1\fi\let\next=\f@gs\fi\next}
%
%
\newwrite\lfile

{\escapechar-1\xdef\pctsign{\string\%}\xdef\leftbracket{\string\{}
\xdef\rightbracket{\string\}}\xdef\numbersign{\string\#}}

\def\writestop{\def\writestoppt{\immediate\write\lfile{\string\pageno%
\the\pageno\string\startrefs\leftbracket\the\refno\rightbracket%
\string\def\string\secsym\leftbracket\secsym\rightbracket%
\string\secno\the\secno\string\meqno\the\meqno}\immediate\closeout\lfile}}
\def\writestoppt{}\def\writedef#1{}

\def\seclab#1{\xdef #1{\the\secno}\writedef{#1\leftbracket#1}\wrlabel{#1=#1}}

\def\subseclab#1{\xdef #1{\secsym\the\subsecno}%
\writedef{#1\leftbracket#1}\wrlabel{#1=#1}}
\def\appsubseclab#1{\xdef #1{\secsym\the\appsubsecno}%
\writedef{#1\leftbracket#1}\wrlabel{#1=#1}}
\def\subsubseclab#1{\xdef #1{\secsym\the\subsecno.\the\subsubsecno}%
\writedef{#1\leftbracket#1}\wrlabel{#1=#1}}
\newwrite\tfile \def\writetoca#1{}
\def\leaderfill{\leaders\hbox to 1em{\hss.\hss}\hfill}
\def\writetoc{\immediate\openout\tfile=toc.tmp
   \def\writetoca##1{{\edef\next{\write\tfile{\noindent ##1
   \string\leaderfill {\noexpand\number\pageno} \par}}\next}}}

%
\catcode`\@=12 
%
%
%
%
%
\def\dbend{{\manual\char127}}
\def\d@nger{\medbreak\begingroup\clubpenalty=10000
    \def\par{\endgraf\endgroup\medbreak} \noindent\hang\hangafter=-2
    \hbox to0pt{\hskip-\hangindent\dbend\hfill}\ninepoint}
\outer\def\danger{\d@nger}

\def\p{\partial}

\def\darr#1{\raise1.5ex\hbox{$\leftrightarrow$}\mkern-16.5mu #1}
\def\half{{\textstyle{1\over2}}} 

%
%
\def\al{\alpha}
\def\be{\beta}
  
\def\de{\delta}  \def\De{\Delta}
\def\ep{\epsilon}

\def\th{\theta}    

\def\ka{\kappa}
\def\la{\lambda} \def\La{\Lambda}
\def\rh{\rho}
  
\def\ta{\tau}
  
\def\ph{\phi}  \def\Ph{\Phi}

\def\om{\omega}  
%
%

%

%
%

 \def\cF{{\cal F}}

\def\cM{{\cal M}}

\def\cO{{\cal O}}
\def\cP{{\cal P}}
\def\cR{{\cal R}}

\def\proof{\noindent {\it Proof:}\ }
\def\Box{\hbox{$\rlap{$\sqcup$}\sqcap$}}

%

%
%
\def\amsyes{y }

\def\answ{y }

\ifx\answ\amsyes
\input amssym.def


\def\CC{{\Bbb C}}
\def\ZZ{{\Bbb Z}}
\def\NN{{\Bbb N}}

\def\bfg{{\frak g}}
\def\bfh{{\frak h}}

\def\hg{{\widehat{\frak g}}}
\def\whg{\hg}
\def\sln{\frak{sl}_N}   \def\hsln{\widehat{\frak{sl}_N}}
\def\sltw{\frak{sl}_2}  \def\hsltw{\widehat{\frak{sl}_2}}

\def\gln{\frak{gl}_n}   
\def\glN{\frak{gl}_N}

\else
\def\ZZ{{Z\!\!\!Z}}              
\def\CC{{I\!\!\!\!C}}
\def\NN{{I\!\!N}}

\def\bfg{{\bf g}}
\def\bfh{{\bf h}}

\def\hg{\hat{\bf g}}

\def\sln{s\ell_n}   \def\hsln{\widehat{s\ell_n}}
\def\sltw{s\ell_2}  \def\hsltw{\widehat{s\ell_2}}

\def\gln{g\ell_n}

\def\cR{{\cal R}} 
\fi
%

%
%

%
%

\newsymbol\ltimes 226E
\newsymbol\rtimes 226F
%
%
%

\def\CMP#1{Comm.\ Math.\ Phys.\ {\bf #1}}

\def\IJMP#1{Int.\ J.\ Mod.\ Phys.\ {\bf #1}}

\def\JPA#1{J.\ Phys.\ {\bf A{#1}}}
\def\JRAM#1{J.\ reine angew.\ Math.\ {\bf {#1}}}

\def\LEM#1{L'Enseignement Math\'ematique {\bf {#1}}}
\def\LMP#1{Lett.\ Math.\ Phys.\ {\bf #1}}

\def\MPL#1{Mod.\ Phys.\ Lett.\ {\bf #1}}
\def\NPB#1{Nucl.\ Phys.\ {\bf B#1}}
\def\PLB#1{Phys.\ Lett.\ {\bf {#1}B}}

\def\PRL#1{Phys.\ Rev.\ Lett.\ {\bf #1}}
\def\PTP#1{Prog.\ Theor.\ Phys.\ Suppl.\ {\bf #1}}
\def\SMD#1{Sov.\ Math.\ Dokl.\ {\bf {#1}}}

%

%
%
\def\SMu{\hbox{\lower 3pt\hbox{ \epsffile{su10.eps}}}}
\def\SMs{\hbox{\lower 3pt\hbox{ \epsffile{ss10.eps}}}}
\def\SMd{\hbox{\lower 3pt\hbox{ \epsffile{sd10.eps}}}}

\def\SMS{\leavevmode\vadjust{\rlap{\smash%
{\line{{\escapechar=` \hfill\rlap{\hskip.3in%
                 \hbox{\lower 2pt\hbox{\epsffile{sd10.eps}}}}}}}}}}
\def\SMH{\leavevmode\vadjust{\rlap{\smash%
{\line{{\escapechar=` \hfill\rlap{\hskip.3in%
                 \hbox{\lower 2pt\hbox{\epsffile{su10.eps}}}}}}}}}}
%
%
\def\LW#1{\lower .5pt \hbox{$\scriptstyle #1$}}
\def\LWr#1{\lower 1.5pt \hbox{$\scriptstyle #1$}}
\def\LWrr#1{\lower 2pt \hbox{$\scriptstyle #1$}}
\def\RSr#1{\raise 1pt \hbox{$\scriptstyle #1$}}

%

\hfuzz=16pt
\nopagenumbers
\pageno=0
%
%
\def\ni{\noindent}

\def\sln{\frak{sl}_n}  \def\hsln{\widehat{\frak{sl}_n}}
\def\qbin#1#2{ \left[ \matrix{ #1 \cr #2 \cr} \right] }

\font\manual=manfnt

%
%

\line{}
\vskip2cm
\centerline{\bf SPINON DECOMPOSITION AND YANGIAN STRUCTURE}\smallskip
\centerline{{\bf OF $\hsln$ MODULES}\footnote{$^\dagger$}{To appear in
``Geometric Analysis and Lie Theory in Mathematics and Physics'', 
Lecture Notes Series of the Australian Mathematical Society.}}
\vskip1cm

\centerline{Peter BOUWKNEGT$\,^{1}$ and Kareljan SCHOUTENS$\,^{2}$}
\bigskip

\centerline{\sl $^1$ Department of Physics and Mathematical Physics}
\centerline{\sl University of Adelaide}
\centerline{\sl Adelaide, SA~5005, AUSTRALIA}
\bigskip

\centerline{\sl $^2$ Institute for Theoretical Physics}
\centerline{\sl University of Amsterdam}
\centerline{\sl Valckenierstraat 65}
\centerline{\sl 1018~XE~~Amsterdam, THE NETHERLANDS}
\medskip
\vskip1.5cm

\centerline{\bf ABSTRACT}\medskip
{\rightskip=1cm 
\leftskip=1cm 
\noindent We review several fermionic-type character formulae for
the characters of the integrable highest weight modules of $\hsln$
at level $\ell=1$, and explain how they arise from a spinon basis 
for these modules.  We also review how the 
Yangian $Y(\sln)$ acts on the integrable $(\hsln)_1$ modules and we
decompose the characters with respect to this $Y(\sln)$ action.}

\vfil

\line{ADP-96-36/M47 \hfil}
\line{ITFA-97-08 \hfil}
\line{{{\tt q-alg/9703021}}\hfil }

\eject


\bref\AL{
I.~Affleck and A.W.W.~Ludwig, {\it Critical theory of overscreened 
Kondo fixed points}, \NPB{360} (1991) 641-696; 
{\it Exact critical theory of the two impurity Kondo model},
\PRL{68} (1992) 1046-1049;
A.W.W.~Ludwig and I.~Affleck, 
{\it Exact conformal field theory results on the multichannel
Kondo effect: Asymptotic three-dimensional space and time dependent 
multipoint and many particle Green's functions},
\NPB{428} (1994) 545-611.}

\bref\AN{
C.~Ahn and S.~Nam, {\it Yangian symmetries in the $SU(N)_1$ WZW model
and the Calogero-Sutherland model}, \PLB{378} (1996) 107-112,
({\tt hep-th/9510242}).}

\bref\An{
G.E.~Andrews, {\it The theory of partitions}, (Addison-Wesley, Reading, 1976).
}

\bref\ANOT{
T.~Arakawa, T.~Nakanishi, K.~Oshima and A.~Tsuchiya,
{\it Spectral decomposition of the path space in solvable lattice model},
{\tt (q-alg/9507025)}.}

\bref\Be{
A.~Berkovich, {\it Fermionic counting of RSOS states and Virasoro 
character formulas for the unitary minimal series $M(\nu,\nu+1)$},
\NPB{431} (1994) 315-348, {\tt (hep-th/9403073)};
A.~Berkovich and B.~McCoy, {\it Continued fractions and fermionic
representations for characters of $M(p,p')$ minimal models}, \LMP{37} 
(1996) 49-66, {\tt (hep-th/9412030)}.}

\bref\Ber{
D.~Bernard, {\it An introduction to Yangian symmetries}, 
({\tt hep-th/9211133}).}

\bref\BGHP{
D.~Bernard, M.~Gaudin, F.D.M.~Haldane and V.~Pasquier,
{\it Yang-Baxter equation in long range interacting systems},
\JPA{26} (1993) 5219-5236, ({\tt hep-th/9301084}).}

\bref\BPS{
D.~Bernard, V.~Pasquier and D.~Serban, {\it Spinons in conformal
field theory}, \NPB{428} (1994) 612-628, ({\tt hep-th/9404050}).}

\bref\BLSa{
P.~Bouwknegt, A.~Ludwig and K.~Schoutens, 
{\it Spinon bases, Yangian symmetry and fermionic representations of 
Virasoro characters in conformal field theory},
\PLB{338} (1994) 448-456, ({\tt hep-th/9406020}).}

\bref\BLSb{
P.~Bouwknegt, A.~Ludwig and K.~Schoutens, 
{\it Spinon basis for higher level $SU(2)$ WZW models}, 
\PLB{359} (1995) 304-312, ({\tt hep-th/9412108}).}

\bref\BS{
P.~Bouwknegt and K.~Schoutens, {\it The $\widehat{SU(n)}_k$ WZW models:
Spinon decomposition and Yangian structure}, 
\NPB{482} (1996) 345-372, ({\tt hep-th/9607064}).}

\bref\CPa{
V.~Chari and A.~Pressley, {\it Yangians and $R$-matrices}, 
\LEM{36} (1990) 267-302.}

\bref\CPb{
V.~Chari and A.~Pressley, {\it Fundamental representations of Yangians 
and singularities of $R$-matrices}, \JRAM{417} (1991) 87-128.}

\bref\CPc{
V.~Chari and A.~Pressley, {\it A guide to Quantum Groups},
(Cambridge University Press, Cambridge, 1994).}
 
\bref\DJKMO{
E.~Date, M.~Jimbo, A.~Kuniba, T.~Miwa and M.~Okado, {\it Maya diagrams
and representation of $\frak{sl}(r,\CC)$}, Adv.\ Stud.\ in Pure Math.\
{\bf 19} (1989) 149-191.}

\bref\Dra{
V.G.~Drinfel'd, {\it Hopf algebras and the quantum Yang-Baxter equation},
\SMD{32} (1985) 254-258.}

\bref\Drb{
V.G.~Drinfel'd, {\it A new realization of Yangians and quantized affine
algebras}, \SMD{36} (1988) 212-216.}

\bref\Drc{
V.G.~Drinfel'd, {\it Quantum groups}, in Proc.\ of the International
Congress of Mathematicians, Berkeley, 1986.}

\bref\FS{
B.L.~Feigin and A.V.~Stoyanovsky, {\it Quasi-particle models for the
representations of Lie algebras and geometry of flag manifolds},
({\tt hep-th/9308079}).}

\bref\FLS{
P.~Fendley, A.W.W.~Ludwig and H.~Saleur, {\it Exact conductance 
through point contacts in the $\nu=1/3$ fractional quantum Hall
effect}, \PRL{74} (1995) 3005, {\tt (cond-mat/9408068)}.}

\bref\FQ{
O.~Foda and Y.-H. Quano, 
{\it Polynomial identities of the Rogers--Ramanujan type},
{\tt (hep-th/9407191)}; {\it Virasoro character identities from the
Andrews--Bailey construction}, {\tt (hep-th/9408086)}.}

\bref\Gea{
G.~Georgiev, {\it Combinatorial constructions of modules for 
infinite-dimensional Lie algebras, I.\ Principal subspace},
({\tt hep-th/9412054}).}

\bref\Geb{
G.~Georgiev, {\it Combinatorial constructions of modules for 
infinite-dimensional Lie algebras, II.\ Parafermionic space}, 
({\tt q-alg/9504024}).}

\bref\HHTBP{
F.D.M.~Haldane, Z.N.C.~Ha, J.C.~Talstra, D.~Bernard and V.~Pasquier,
{\it Yangian symmetry of integrable quantum chains with long-range 
interactions and a new description of states in conformal field 
theory}, \PRL{69} (1992) 2021-2025.}

\bref\Hik{
K.~Hikami, {\it Yangian symmetry and Virasoro character in a lattice
spin system with long-range interactions},
\NPB{441} (1995) 530-548.}

\bref\JKKMP{
M.~Jimbo, R.~Kedem, H.~Konno, T.~Miwa and J.-U.~Petersen, 
{\it Level-$0$ structure of level-$1$ $U_q(\hsltw)$-modules
and Macdonald polynomials}, \JPA{28} (1995) 5589-5606,
({\tt q-alg/9506016}).}

\bref\KKMM{
R.~Kedem, T.~Klassen, B.~McCoy and E.~Melzer, {\it Fermionic 
quasiparticle representations for characters of $G^{(1)}\times 
G^{(1)} / G^{(2)}$},
\PLB{304} (1993) 263-270, {\tt (hep-th/9211102)}; {\it Fermionic sum 
representations for conformal field theory characters},
\PLB{307} (1993) 68-76, {\tt (hep-th/9301046)};
S.~Dasmahapatra, R.~Kedem, T.~Klassen, B.~McCoy and E.~Melzer, 
{\it Quasiparticles, conformal field theory and $q$ series}, 
\IJMP{B7} (1993) 3617-3648, {\tt (hep-th/9303013)};
R.~Kedem, B.~McCoy and E.~Melzer, 
{\it The sums of Rogers, Schur and Ramanujan and the Bose-Fermi
correspondence in $1+1$-dimensional quantum field theory},
in ``Recent progress in Statistical Mechanics and Quantum Field Theory'',
p.195-219, (World Scientific, Singapore, 1995),
{\tt (hep-th/9304056)};
E.~Melzer, {\it The many faces of a character}, \LMP{31} (1994) 233-246,
{\tt (hep-th/9312043)}.}

\bref\Kir{
A.N.~Kirillov, {\it Dilogarithm identities}, \PTP{118} (1995) 61-142,
({\tt hep-th/9408113}).}

\bref\KKN{
A.N.~Kirillov, A.~Kuniba and T.~Nakanishi, {\it Skew Young diagram method 
in spectral decomposition of integrable lattice models}, 
({\tt q-alg/9607027}).}

\bref\KNS{
A.~Kuniba, T.~Nakanishi and J.~Suzuki,
{\it Characters in conformal field theories from thermodynamic Bethe 
Ansatz}, \MPL{A8} (1993) 1649-1660, {\tt (hep-th/9301018)}.}

\bref\MD{
I.G.~Macdonald, {\it Symmetric functions and Hall polynomials}, 2nd
edition, (Oxford University Press, New York, 1995).}

\bref\NYa{
A.~Nakayashiki and Y.~Yamada, {\it Crystalizing the spinon basis},
\CMP{178} (1996) 179-200, ({\tt hep-th/9504052}).}

\bref\NYb{
A.~Nakayashiki and Y.~Yamada, {\it Crystalline spinon basis for RSOS
models}, \IJMP{A11} (1996) 395-408, ({\tt hep-th/9505083}).}


\bref\NT{
M.~Nazarov and V.~Tarasov, {\it Representations of Yangians with
Gel'fand-Zetlin basis}, ({\tt q-alg/9502008}).}

\bref\Ol{
G.I.~Olshanskii, {\it Representations of infinite-dimensional classical
groups, limits of enveloping algebrass, and Yangians}, in ``Topics
in Representation Theory,'' Adv.\ Sov. Math. {\bf 2}, 1-66 (Amer.\ 
Math.\ Soc., Providence, 1991).}

\bref\Sc{
K.~Schoutens, {\it Yangian symmetry in conformal field theory}, 
\PLB{331} (1994) 335-341, ({\tt hep-th/9401154}).}

\bref\TUa{
K.~Takemura and D.~Uglov, {\it The orthogonal eigenbasis and norms 
of eigenvectors in the spin Calogero-Sutherland model}, 
({\tt solv-int/9611006}).}


\bref\Ug{
D.~Uglov, {\it Yangian Gelfand-Zetlin basis, $\glN$-Jack polynomials
and computation of dynamical correlation functions in the 
spin Calogero-Sutherland model},
({\tt hep-th/9702020}).}

%
\baselineskip=1.5\baselineskip

\footline{\hss \tenrm -- \folio\ -- \hss}


\newsec{Introduction}

Conventionally, the Hilbert space 
of (rational) two dimensional conformal field theories (RCFT)
is described in terms of a chiral algebra that acts on a finite set of
fields that are primary with respect to this chiral algebra.
This procedure leads to so-called Verma-module bases of RCFT's,
and gives rise to `bosonic type' formulas for the characters
(i.e.\ torus partition functions) of such conformal field theories,
reflecting this particular choice for the basis of the Hilbert space.
On the other hand, in a conformal field theory (CFT) there
are many bases of Hilbert space, none of which is a priori
preferred. Which of those different bases is most
useful depends on the question that is being asked (see e.g.\
[\AL,\FLS]). This is quite different in a massive
field theory, where a basis of massive particles is 
naturally distinguished. In general, considering massive
perturbations of CFT's, and their corresponding particle
bases, one may generate, at least conceptually, ``quasi-particle''
bases of CFT's, by letting the mass tend to zero. This provides
a way to define the notion of a massless ``quasi-particle'',
which has proven to be very useful 
in cases where the massive perturbation, used to define
such a basis, is Yang-Baxter integrable.
The so-generated bases inherit a particle (``Fock-space'') like 
structure, which is not manifest
in the Verma-module type basis of the same CFT.
Each basis of the Hilbert space
of a given CFT gives rise to a particular way of writing
the  partition function on a torus.
The equality of the different ways to write 
the torus partition function, in different bases, gives rise
to remarkable identities. Moreover, bases of massless
quasiparticles appear to be deeply related to Yangian and affine
quantum symmetries, that are often present even when
conformal invariance is broken. Therefore, we expect
that a better understanding of such bases of CFT's will also
provide valuable insights into perturbed CFT's.
Analysis of the thermodynamic Bethe Ansatz, arising
from (Yang-Baxter) integrable perturbations of CFT's, for a variety
of models, has  led to a wealth of conjectures for so-called
quasi-particle (or fermionic) type characters for  
conformal field theories (see, in particular [\KNS,\KKMM]),
some of which have been proven through 
$q$-analysis (see e.g.\ [\Be,\Kir,\FQ]).
Most of these results still lack an
interpretation (and/or proof) in terms of a corresponding 
structure of the Hilbert space of a 
conformal field theory, i.e.\
a characterization  and/or construction of the corresponding ``quasi-particle''
basis of the Hilbert space.     

A particularly interesting model illustrating the issues above 
is the so-called $\sln$ Haldane-Shastry long-range spin chain, which is 
integrable and has Yangian symmetry (even for finite chains).
Its low energy sector is identical to
a well-known conformal field 
theory, namely the $SU(n)$ level-$1$ Wess-Zumino-Witten model [\HHTBP,\BGHP]. 
While the description of the Hilbert space 
of the $SU(n)$ level-$1$ WZW model in terms of its chiral 
algebra, i.e.\ $\hsln$, is complicated due to the existence of
null vectors, it was found that 
the Hilbert space has a very simple structure [\BPS,\BLSa,\BS], originating
from the underlying Yangian symmetry: it
may be viewed as a ``Fock space'' of massless `spinon particles', which
satisfy generalized commutation relations [\BLSa] (and no other relations
for $\hsltw$).
For generalizations of these results to higher level $\ell>1$ 
we refer to [\BLSb,\NYa,\NYb,\ANOT].

In this paper we review these new exciting developments from a mathematical
point of view.  

Let us illustrate the 
various issues in the case of the affine Lie algebra 
$\hsltw$.  The bosonic character formula for the level-$1$ integrable
highest weight modules $L(\widehat{\La}_k),\, k=0,1$, reads
(see App.\ B)
\eqn\eqAAc{
{\rm ch}_{L(\widehat{\La}_k)}(z;q) ~=~ \sum_{n\in\ZZ} \,
  {q^{(n+ {1\over2}k)^2} \over (q)_\infty} z^{n+{k\over2}}\,.
}
There are two, not obviously related, fermionic character formulae
\eqn\eqAAa{
{\rm ch}_{L(\widehat{\La}_k)}(z;q) ~=~ q^{ {1\over4} k^2 }\, 
  \sum_{m_1,m_2\geq 0} { q^{m_1^2 - m_1m_2 + m_2^2 + k(m_1-m_2)} 
  \over (q)_{m_1} (q)_{m_2} } z^{m_1 - m_2 + {1\over2} k}\,,
}
and 
\eqn\eqAAb{ \eqalign{
{\rm ch}_{L(\widehat{\La}_k)}(z;q) & ~=~
  \sum_{m_1,m_2\geq 0\atop m_1 + m_2 \equiv k\, {\rm mod}\,2} 
  { q^{{1\over 4} (m_1 + m_2)^2} \over (q)_{m_1}
  (q)_{m_2} } z^{{1\over2}(m_1 - m_2)}\cr
& ~=~
  \sum_{m_1,m_2\geq 0\atop m_1 + m_2 \equiv k\, {\rm mod}\,2} 
  \sum_{m\geq0} \, (-1)^m\, {q^{{1\over2} m(m-1)} \over (q)_m}\,
  { q^{{1\over 4} (m_1 - m_2)^2} \over (q)_{m_1-m}
  (q)_{m_2-m} } z^{{1\over2}(m_1 - m_2)}\,,\cr}
}
which are easily shown to be equivalent to \eqAAc\ by using the 
Durfee square (Lemma D.2).  The second equality in \eqAAb\ follows
from Lemma D.3 (i).  [As compared to the formulae in Appendix B we
have defined $z = {x_1\over x_2}$.]

While \eqAAa\ clearly has an interpretation in terms of two quasi-particles
associated to the roots of $\sltw$, equation \eqAAb\ arises from two
quasi-particles associated to the weights of the fundamental 
two-dimensional representation of $\sltw$ (i.e., the spinor representation of
$\frak{so}_3$, hence the name spinon).  Moreover, since the quasi-particle
operator intertwines between the modules $L(\widehat{\La}_0)$ and 
$L(\widehat{\La}_1)$,  applying an even or odd number of creation operators 
to the vacuum results in a vector in $L(\widehat{\La}_0)$ or 
$L(\widehat{\La}_1)$, respectively.
In this review we will only
consider the generalization of \eqAAb\ to $\sln$.  For the 
generalization of \eqAAa\ to $\sln$, see [\FS,\Gea,\Geb].
The spinon operators are mutually nonlocal.  As a consequence they 
satisfy generalized commutation relations which leads to 
fractional statistics.  Nevertheless, the generalized commutation relations
are sufficiently powerful to determine a linearly independent set of
basis vectors for the spinon Fock space.  Calculating the character
of the $(\hsltw)_1$ integrable highest weight modules using the spinon 
basis immediately leads to \eqAAb.

 The integrable highest weight modules $L(\widehat{\La}_k),\, k=0,1$,
of $(\hsltw)_1$
admit a (semi-simple) action of a certain quantum group, the so-called
Yangian $Y(\sltw)$.  This action is not only most naturally described 
on the spinon basis but is, in fact, intimately related to the existence
of a spinon basis.  The $1$-spinon states constitute the 
(two-dimensional) evaluation representation of $Y(\sltw)$ with the spinon
mode index playing the role of the evaluation parameter, while the action 
on multi-spinon states corresponds to a (non-standard) co-multiplication
for $Y(\sltw)$.

The modules $L(\widehat{\La}_k),\, k=0,1$, can be decomposed with respect 
to this $Y(\sltw)$ action.  For $\sltw$ this leads to the characterformula
\eqn\eqAAf{
{\rm ch}_{L(\widehat{\La}_k)}(z;q)  ~=~
  \sum_{N\geq0\atop N\equiv k\,{\rm mod\,}2}\ 
  \sum_{\la\in\cP\atop  l(\la)\leq N}\ q^{ {N^2\over4} + |\la| } \ 
  {\rm ch}_{\la}^{Y(\sltw)} (z) \,,
}
where, for each spinon number $N$, the sum is over all partitions
(Young diagrams)
$\la$ of length $l(\la)\leq N$ and ${\rm ch}_{\la}^{Y(\sltw)} (z)$
denotes the character of an irreducible finite-dimensional $Y(\sltw)$
module $L_\la$ labeled by $\la$ (see Section 3 for more details).
The generalization of this result to $(\hsln)_1$ turns out to be slightly 
more complicated; in that case the Yangian modules that occur are 
most naturally parametrized in terms of a particular kind of skew Young
diagrams, the so-called border strips.

This survey is organized as follows.  In Section 2 we prove various 
fermionic-type formulae for the characters of the integrable modules 
of $\hsln$ at level $\ell=1$, and explain their origin in terms of a 
spinon basis for these modules.  In Section 3 we explain how the 
Yangian $Y(\sln)$ acts on the integrable $(\hsln)_1$ modules and we
decompose the characters with respect to this $Y(\sln)$ action.
In order not to deter the reader from the main line of thought we
have omitted most of the proofs (which can either be found in the 
literature or are straightforward) and 
deferred most of the mathematical prerequisites to appendices.
Appendix A is a brief introduction
to partitions and symmetric functions. In particular we  introduce
(skew) Schur functions.  
Appendix B briefly introduces the Lie algebra $\sln$, its (untwisted) affine
extension $\hsln$ and some of their modules.  
In Appendix C we  
define the Yangian $Y(\sln)$ and discuss its finite-dimensional 
irreducible representations and, finally, in Appendix D we discuss some
$q$-identities necessary to establish the character identities of Section 2.

\newsec{Spinon decomposition of $\hsln$ modules}

\subsec{$N$-spinon cuts of the $(\hsln)_1$ characters}

In this section we prove two fermionic-type character formulae
for the characters of the integrable highest weight modules 
$L(\widehat{\La}_k)$, $k=0,\ldots,n-1$, of the affine Lie algebra
$\hsln$ at level $\ell=1$.  In the next section we will then
argue that these formulae find their origin in the existence 
of a certain basis, the so-called spinon basis, for these modules,
where we identify the summation variable $N$ as the spinon number.
We refer to Appendix B for definitions and notations regarding
$\sln$ and its affinization $\hsln$. 
\thm\thBAa
\proclaim Theorem \thBAa.  For $\widehat{\La}\in P_+^{\,1}$ and 
$\La- \la\in Q$
we have the following expressions for the $(\hsln)_1$ string functions
$c_\la^\La(q)$ defined by
\eqn\eqBAaaa{
{\rm ch}_{L(\widehat{\La})} ~=~ \sum_\la \, c_\la^\La(q) e^\la\,.
}
\eqn\eqBAaa{
c^\La_\la(q) ~=~ {q^{{1\over2} |\la|^2 } \over (q)_\infty^{n-1} } ~=~ 
  \sum_{N\geq0} \, 
  c^{\La, N}_\la(q)\,,
}
where 
\eqn\eqBAab{ \eqalign{
c_\la^{\La,N} & ~=~ q^{{1\over2} |\la|^2} 
  \sum_{m_1,\ldots,m_{n-2}} \ 
  { q^{A_1m_1 + (A_2-m_1)m_2 +\ldots + (A_{n-1} - (m_1+\ldots+m_{n-2}))
    (A_{n} - (m_1+\ldots+m_{n-2}))} \over (q)_{A_1} (q)_{A_2-m_1} \ldots
    (q)_{A_{n-1} - (m_1+\ldots+m_{n-2})} (q)_{A_{n} - (m_1+\ldots+m_{n-2})}}\cr
  &\quad \quad \times \, {1\over \prod_i (q)_{m_i}} \cr
  & ~=~ q^{{1\over2} |\la|^2} 
  \sum_m\ (-1)^m\, { q^{{1\over2}m(m-1) } \over (q)_m} \,
  {1\over \prod_i (q)_{A_i -m} }\,,\cr}
}
for those $N$ such that $\widehat{\La} = \widehat{\La}_{N\, {\rm mod}\,n}$
and zero otherwise.
We have defined 
\eqn\eqBAac{
A_i ~=~ {N\over n} + (\la,\ep_i) \,.
}
\par

\proof  
Consider the first expression for $c_\la^{\La,N}$.  Changing 
$N\to N+ (m_1+\ldots +m_{n-2})n$ we have 
$$\eqalign{
\sum_{N\geq0} \,c^{\La, N}_\la(q) ~=~ & q^{{1\over2} |\la|^2} 
  \sum_{N}\sum_{m_1,\ldots,m_{n-2}} \ 
  {q^{(A_1+m_1+\ldots+m_{n-2})m_1 + (A_2+m_2+\ldots+m_{n-2})m_2
  + \ldots + A_{n-1}A_n } \over (q)_{A_1+m_1+\ldots+m_{n-2}}
  (q)_{A_2+m_2+\ldots+m_{n-2}}\ldots (q)_{A_{n-1}} (q)_{A_n}}\cr
& \times \, {1\over \prod_i (q)_{m_i} } \cr}
$$
Then apply the Durfee square (Lemma D.2) to, respectively, $m_1, m_2, \ldots
m_{n-2}$ and $N$ to obtain \eqBAaa.
To prove the second expression in \eqBAab\ we 
successively apply Lemma D.3 (ii)
to the right hand side, i.e.\ 
$$ \eqalign{
\sum_{m} &  (-1)^m\, { q^{{1\over2}m(m-1) } \over (q)_m} \,
  {1\over \prod_i (q)_{A_i -m} } \cr
& ~=~ \sum_{m, m_{n-2}}  (-1)^m\, { q^{{1\over2}m(m-1) } \over (q)_m} \,
  {q^{(A_{n-1}-m_{n-2}-m)(A_{n}-m_{n-2}-m)}\over
  (q)_{A_1-m} \ldots (q)_{A_{n-2}-m} (q)_{A_{n-1}-m_{n-2}-m} 
  (q)_{A_{n}-m_{n-2}-m} } \cr
& ~=~ ~~\ldots\ldots \cr
& ~=~ \sum_{m,m_1,\ldots,m_{n-2}} \ (-1)^m\, { q^{{1\over2}m(m-1) } 
  \over (q)_m} \, \cr
& \quad ~\times~ 
  {q^{(A_2-m_1-m)(m_2-m_1) + (A_3-m_2-m)(m_3-m_2) + \ldots
  + (A_{n-1}-m_{n-2}-m)(A_{n}-m_{n-2}-m)} \over 
  (q)_{A_1-m} (q)_{A_2-m_1-m}(q)_{A_3-m_2-m} \ldots (q)_{A_{n}-m_{n-2}-m}} \cr
& \quad
  ~\times~ {1\over (q)_{m_1} (q)_{m_2-m_1} \ldots (q)_{m_{n-2}-m_{n-3}}} \cr}
$$
Now change summation variables
$$ \eqalign{
m_1 & ~\to~ m_1 - m \cr
m_2 & ~\to~ m_2 + m_1 - m \cr
\vdots & ~~~~~~~~~~~\vdots \cr
m_{n-2} & ~\to~ m_{n-2} + m_{n-3} + \ldots  + m_1 - m \cr}
$$
and use Lemma D.3 (i) in the form
$$
\sum_{m} \ (-1)^m\, { q^{{1\over2}m(m-1) } 
  \over (q)_m} { 1 \over (q)_{A_1-m} (q)_{m_1 -m} } ~=~ { q^{A_1 m_1} \over
  (q)_{A_1} (q)_{m_1} }\,,
$$
to obtain \eqBAab.\Box

Note that, for $n=2$, the formulae \eqBAab\ lead to those of
\eqAAb.


\subsec{Generalized commutation relations and spinon basis}

The character formulae of Theorem \thBAa\ find their origin in the 
existence of a specific basis, the so-called spinon basis, for 
the integrable highest weight modules of $\hsln$ at level $\ell=1$.
We will refer to the intertwiners
\eqn\eqBBaa{
\ph^i \left( \matrix{ \La_1 \cr \La_{k+1} \ \La_k \cr} \right) (z)
 \ : \ L(\widehat{\La}_k)~\to~ L(\widehat{\La}_{k+1})\,,\qquad (i=1,\ldots,n)
\,,
}
transforming in the fundamental $\sln$ representation $L(\La_1)$
(see Appendix B for a more precise statement), as the spinon operators
of $\hsln$.  (Here the subscripts on $\widehat{\La}$ are taken modulo
$n$.) 
Their mode expansion is given by
\eqn\eqBBab{
\ph^i \left( \matrix{ \La_1 \cr \La_{k+1} \ \La_k \cr} \right) (z)
 ~=~ \sum_{m\in\ZZ}\, 
  \ph^i \left( \matrix{ \La_1 \cr \La_{k+1} \ \La_k \cr} \right)_
  {- {n-2k-1\over 2n} - m} z^{m - {k\over n}}\,.
}
Multiple application of the spinon creation operators to the 
highest weight state $|0\rangle$ of $L(\widehat{\La}_0)$ creates a spinon Fock
space.  By definition, a state with $N$ spinon excitations is a vector
in $L(\widehat{\La}_{N\,{\rm mod}\, n})$.  Hence, it is unambiguous on
which space $L(\widehat{\La}_k)$ the spinon operator is acting and 
henceforth we will omit the fusion channel and simply write $\ph^i(z)$.

A general $N$-spinon state can be written as a linear combination of vectors
\eqn\eqBBac{
|n_1,\ldots, n_N\rangle ~\equiv~ \ph^{i_N}_{- {n-(2N-1)\over 2n} -n_N} \ldots
  \ph^{i_3}_{- {n-5\over2n} -n_3} 
  \ph^{i_2}_{- {n-3\over2n} -n_2}
  \ph^{i_1}_{- {n-1\over2n} -n_1} |0\rangle\,.
}
The collection of all the vectors \eqBBac\
forms an (overcomplete) basis for the $(\hsln)_1$
module $\oplus_{k=0}^{n-1} \, L(\widehat{\La}_k)$.  Note that the energy
of the state \eqBBac\ is given by
\eqn\eqBBaf{
L_0 |n_1,\ldots, n_N\rangle ~=~ \left( {N(n-N)\over 2n} + \sum_{i=1}^N n_i
 \right) |n_1,\ldots, n_N\rangle\,.
}

Due to the mutual nonlocality of the spinon operators, the modes
\eqBBab\ will not satisfy simple commutation relations.  Instead, 
the modes will satisfy so-called generalized commutation relations.
Nevertheless, these  generalized commutation relations imply that 
not all vectors \eqBBac\ are linearly independent.  In particular,
they allow us to choose an ordering $0\leq n_1\leq n_2\leq \ldots
\leq n_N$.  Let us consider the abstract algebra generated by the spinon 
modes modulo the generalized commutation relations and the corresponding
module, referred to a the spinon Fock space $\cF$, spanned by the
vectors \eqBBac.  The problem is to determine a linearly independent 
set of basis vectors for $\cF$.

Let us illustrate this explicitly for $\hsltw$.
In this case the generalized commutation relations for the spinon fields
are (see [\BLSa])
\eqn\eqBBad{
\sum_{l\geq0} C^{(-{1\over2})}_l \left(
  \ph^i_{-p - {k+1\over2} - l + {3\over4}} \ph^j_{-q - {k\over2} + l 
   + {3\over4}} - \left( \matrix{ i\leftrightarrow j \cr p\leftrightarrow
  q\cr} \right)  \right) ~=~ (-1)^k\ep^{ij} \de_{p+q+k-1}\,,
}
where $k=0,1$, depending on whether the left hand side acts on
a vector in $L(\widehat{\La}_k)$, $k=0,1$, and the coefficients $C^{(\al)}_l$
are defined by the expansion
\eqn\eqBBae{
(1-x)^\al ~=~ \sum_{l\geq 0} C_l^{(\al)} x^l\,.
}
By induction we find the following refinement of the basis \eqBBac\
\thm\thBBa
\proclaim Theorem \thBBa\ [\BLSa].  The following (linearly 
independent) set of vectors 
form a basis for the spinon Fock space $\cF$
\eqn\eqBBag{
\ph^2_{- {2(M_1+M_2)-1\over 4} -m_{M_2}} \ldots
\ph^2_{- {2(M_1+1)-1\over 4} -m_{1}}
\ph^1_{- {2(M_1-1)-1\over 4} -n_{M_1}}\ldots
\ph^1_{- {1\over4} - n_1 } |0\rangle \,.
}
The $L_0$ eigenvalue of the state \eqBBag\ is given by
\eqn\eqBBah{
{(M_1+M_2)^2\over4} + \sum_{i=1}^{M_2} m_i + \sum_{i=1}^{M_1} n_i\,.
}
\par

A priori, the spinon Fock space $\cF$ could be bigger than 
$L(\widehat{\La}_0) \oplus L(\widehat{\La}_1)$.  In this case,
calculating the character of $\cF$ (using \eqBBah) 
immediately leads to \eqBAab\ (see also \eqAAb), proving
$\cF \cong L(\widehat{\La}_0) \oplus L(\widehat{\La}_1)$. 
For $\hsln$, $n>2$, however, there might be additional relations 
between the spinon operators beyond those coming from the 
generalized commutation relations, i.e.\ there might be  null vectors
in the Fock module $\cF$.  The details remain to be worked out.
In any case, all relations are encoded in
the character formulae \eqBAab.

\newsec{Yangian structure of $\hsln$ modules}

\subsec{Action of $Y(\sln)$ on integrable $(\hsln)_1$ modules}

Let $x^a$ be an orthonormal basis of $\sln$ and let $f^{abc}$ and
$d^{abc}$ be, respectively, the structure constants and $3$-index 
$d$-symbol of $\sln$ normalized such that 
\eqn\eqCAac{ \eqalign{
f^{abc} f^{dbc} & ~=~ -2n\, \de^{ad}\,,\cr
d^{abc} d^{dbc} & ~=~ {2 (n^2-4) \over n } \, \de^{ad} \,,\cr}
}
i.e., in the fundamental representation $(L(\La_1),\rh)$, we have 
\eqn\eqCAad{
t^a t^b ~=~ \textstyle{{1\over n}} \de^{ab} + \half f^{abc} t^c 
  +\half d^{abc} t^c \,,
}
where $t^a \equiv \rh(x^a)$.

The following theorem was suggested by taking the infinite chain limit 
of the Yangian generators in the $\sln$ Haldane-Shastry spin chain
\thm\thCAa
\proclaim Theorem \thCAa\ [\Sc].  The following formulae define an
(semi-simple) action of $Y(\sln)$, as defined in Definition C.5, 
on the integrable modules of $\hsln$ at level $\ell=1$
\eqn\eqCAaa{ \eqalign{
x^a & ~=~ x_0^a \,,\cr
J(x^a) & ~=~ \half f^{abc} \ \sum_{m>0}\, (x_{-m}^b x_m^c ) - 
  {n\over 2(n+2)} W_0^a   \,.\cr}
}
where 
\eqn\eqCAab{
W_0^a ~=~ \half d^{abc} \ \sum_{m\in\ZZ} \, :x_{-m}^b x_m^c:\,,
}
\par

The proof of Theorem \thCAa\ is straightforward, albeit cumbersome,
and amounts to checking the defining relations (Y1)--(Y3)
of Definition C.5.
 
It is useful to determine the action of the Yangian generators
\eqCAaa\ on the spinon basis \eqBBac.  
To this end, introduce (formal) generating
series for the spinon basis,
\eqn\eqCAae{
\Ph^{i_N,\ldots,i_1}(z_N,\ldots,z_1) ~\equiv~
  \ph^{i_N}(z_N)\ldots \ph^{i_2}(z_2)\ph^{i_1}(z_1)|0\rangle\,.
}
By using the commutation relations between $x_m^a$ and the spinon 
operators (see eqn.\ (B.24)) as well as the null vector structure of
the level-$1$ integrable highest weight modules, one arrives at
\thm\thCAb
\proclaim Theorem \thCAb\ [\BPS,\AN].  The action of the Yangian
generators \eqCAaa\ on the spinon basis is given (up to an automorphism 
of $Y(\sln)$) by
\eqn\eqCAaf{\eqalign{
x^a \Ph(z_N,\ldots,z_1) & ~=~ \sum_i\,  t^a_i \,\Ph(z_N,\ldots,z_1)\,,\cr
J(x^a) \Ph(z_N,\ldots,z_1) & ~=~ \left( -n \sum_i D_i t^a_i + 
  \half f^{abc} \sum_{i\neq j} \th_{ij} t_i^b t_j^c \right) 
\Ph(z_N,\ldots,z_1)\,,\cr}
}
where $D_i = z_i \p_{z_i}$ and $\th_{ij} = {z_i \over z_i - z_j}$,
and $t^a_i$ denotes the action of $t^a$ on the $i$-th entry of $\Ph$.\par

The differential operators above are precisely the Yangian generators
of yet another well-known model with Yangian symmetry, namely the
$\sln$-spin generalization of the Calogero-Sutherland model
(at coupling constant $\be=-\half$) [\BGHP].


\subsec{Decomposition of characters under $Y(\sln)$}

After having established that there exists a semi-simple action of
the Yangian $Y(\sln)$ on the integrable highest weight modules of 
$\hsln$ at level $\ell=1$, it is natural to ask how these modules
decompose under $Y(\sln)$.  This is settled by the following theorem
\thm\thCBa
\proclaim Theorem \thCBa\ [\BS,\KKN].  The character of the level-$1$
irreducible module $L(\widehat{\La}_k)$ of $\hsln$ decomposes under the
action of $Y(\sln)$ as
\eqn\eqCBa{
{\rm ch}_{L(\widehat{\La}_k)}(x;q) ~=~  
\sum_{\ka\in {\rm BS}_n \atop |\ka|\equiv k {\rm mod}\, n} \
q^{ E(\ka) } s_\ka(x) \,,
}
where $\ka$ runs over the set of border strips of rank $n$ 
and parametrizes a set of irreducible (tame) $Y(\sln)$ modules with
character $s_\ka(x)$ (Appendix C), and  
\eqn\eqCBb{\eqalign{
E(\ka) & ~=~ { (n-1)\over 2n} |\ka|^2 + \sum_{i=1}^r (i-r)a_i \cr
& ~=~ { |\ka| (n-|\ka|) \over 2n} + \sum_{i=1}^s (s-i)b_i\,,\cr}
}
if $\ka = \langle a_1,\ldots, a_r\rangle = [b_1,\ldots,b_s]$ with $b_s<n$.
\par

\proof To prove the character equality observe that, in terms of 
semi-infinite rapidity sequences (see Appendix A), we can write
\eqn\eqCBba{
E(\ka) ~=~ \De_k - \sum_i (m_i - m_i^{(k)}) \,.
}
Thus the right hand side of \eqCBa\ can be written as 
\eqn\eqCBbc{ 
q^{\De_k}  \sum_{\underline{m}\in\cR_{n}^{\infty/2}}\
  q^{- \sum_i (m_i - m_i^{(k)}) } s_{\underline{m}}(x)\,.
}
which can be shown to equal ${\rm ch}_{L(\widehat{\La}_k)}(x;q)$
either by invoking the path description of $L(\widehat{\La}_k)$
[\DJKMO] (see also [\KKN]) or, as done in [\BS], by `truncating'
to finite rapidity sequences and using a similar result for
the Haldane-Shastry spin chain.  Specifically, \eqCBbc\ equals
\eqn\eqCBbb{ \eqalign{
& \lim_{L\to\infty} q^{\De(\La_{L\,{\rm mod}\,n})}  
  \sum_{\underline{m}\in\cR_{n}^{L}}\
  q^{- \sum_i (m_i - m_i^{(k)}) } s_{\underline{m}}(x) \cr
& ~\equiv~ \lim_{L\to\infty} {\rm ch}_{{\rm HS}}^{(L)}(x;q)\,.\cr}
}
where we have introduced the character of the $\sln$ Haldane-Shastry
spin chain on a chain of length $L$ [\Hik,\BS].
The proof is then completed by taking the $L\to\infty$ limit of 
the following 
\thm\thCBb
\proclaim Lemma \thCBb\ [\BS].  
\eqn\eqCBc{ 
{\rm ch}_{{\rm HS}}^{(L)}(x;q) ~=~ q^{ {n-1\over 2n} L^2 } \,  
  H_L^{(n)}(x;q^{-1})\,,
}
with $H_L^{(n)}(x;q)$ as in (D.5).
\par

In Appendix C (i.e., Theorem C.13) it is shown that the skew Schur 
function $s_\ka(x)$ corresponding to some skew Young diagram
$\ka=\la/\mu$ can be identified with the $\sln$ character of
an irreducible (finite-dimensional) $Y(\sln)$ module.  To prove 
the remainder of Theorem \thCBa, one thus has to construct a $Y(\sln)$ highest
weight vector in $L(\widehat{\La}_k)$ for each $\ka\in{\rm BS}_n$.
This is most easily done in the spinon basis.  For general $n$ the
result can be inferred from [\TUa,\Ug].  Here, for simplicity,
we only consider $\sltw$ [\BPS,\JKKMP].

We quote the following result for $\hsltw$
\thm\thCBd
\proclaim Theorem \thCBd\ [\BPS,\JKKMP].  In the $N$-spinon sector 
of the $(\hsltw)_1$ module $L(\widehat{\La}_0) \oplus L(\widehat{\La}_1)$
there exists a $Y(\sltw)$ highest weight vector $\om_{\la,N}$
for every partition $\la$ with $l(\la)\leq N$.  Furthermore, 
\eqn\eqCBaa{
L_0\, \om_{\la,N} ~=~ \left( |\la| + {N^2\over4}\right)\, \om_{\la,N} \,,
}
and the irreducible $Y(\sltw)$ module $V_{\la,N}$ generated from $\om_{\la,N}$ 
has $\sltw$ character
\eqn\eqCBab{
{\rm ch}_{V_{\la,N}}(x) ~=~ \prod_{i\geq0} \ h_{m_i^\la}(x)\,,
}
where $m_i^\la = \#\{ j: \la_j=i\}$ for $i\geq1$ and $m_0^\la\equiv
N - \sum_{i\geq1} m_i^\la$.\par

As we have seen in Theorem \thCAb, on the spinon basis the Yangian
action is described by the differential operators that generate 
the Yangian symmetry in the spin Calogero-Sutherland model.  For
$\sltw$, it is well-known (see e.g.\ [\BGHP]) that the corresponding
Yangian highest weight vectors are given by Jack polynomials
associated to a partition $\la$.  It 
remains to determine which Jack polynomials can occur in the $N$-spinon
sector (this leads to the restrictions on $\la$) and what the irreducible
$Y(\sltw)$ module is that is generated from this highest weight vector.
This last fact is accomplished by determining the Drinfel'd polynomial 
associated to the highest weight vector (see Appendix C).

To make the comparison of Theorem \thCBd\ with Theorem \thCBa, we have
to associate a skew Young diagram $\ka$ to each partition $\la$ above.
Let $r\in\NN$ be such that $m_i^\la=0$ for $i\geq r$ whilst $m_{r-1}^\la
\neq0$.  Then define the skew Young diagram $\ka =\langle a_1,\ldots,
a_r\rangle$ by
\eqn\eqCBac{
a_i ~=~ \cases{ m_{i-1}^\la +1 & for $i\in \{1,r\}\,,$\cr
                m_{i-1}^\la +2  & for $2\leq i\leq r-1\,.$\cr}
}
Then,
\eqn\eqCBad{\eqalign{
N & ~=~ \sum_{i\geq0} \ m_i^\la ~=~ |\ka| - 2 (r-1)\,,\cr
|\la| & ~=~ \sum_{i\geq1} \ im_i^\la ~=~ \sum_i i a_i - |\ka| - (r-1)^2\,,\cr}
}
such that
\eqn\eqCBae{
|\la| + {N^2\over4} ~=~ {|\ka|^2\over4} + \sum_i \,(i-r)a_i\,,
}
in accordance with \eqCBb.  The equality of the resulting $Y(\sltw)$ 
characters (see \eqCBab) follows from (A.14).\Box

Unfortunately, this simple description of the Yangian highest weight 
vectors in terms of a single partition $\la$ does not
generalize to $\sln$.  
To illuminate the connection of Theorem \thCBa\ to the spinon basis  
for general $n$, consider a generic $M$-spinon state
\eqn\eqCBaf{
|n_1,\ldots, n_M\rangle ~=~ \ph^{i_M}_{- {n-(2M-1)\over 2n} -n_M} \ldots
  \ph^{i_3}_{- {n-5\over2n} -n_3} 
  \ph^{i_2}_{- {n-3\over2n} -n_2}
  \ph^{i_1}_{- {n-1\over2n} -n_1} |0\rangle\,,
}
with $0\leq n_1\leq n_2\leq \ldots \leq n_M$.
Clearly, the energy of such a state is given by (cf.\ \eqBBaf)
\eqn\eqCBag{
L_0 |n_1,\ldots, n_M\rangle ~=~ \left( {M(n-M)\over 2n} + \sum_{i=1}^M n_i
 \right) |n_1,\ldots, n_M\rangle\,.
}
The similarity with \eqCBb\ is striking.  It suggests that to each
sequence $0\leq n_1\leq n_2\leq \ldots \leq n_M$ can be associated a 
skew Young diagram $\ka=[b_1,\ldots,b_s]$ constructed as follows: 
Draw a square, call this
the $1$-st square. Now, by induction, draw the $i$-th square on top 
(resp., to the right) of the $i-1$-th square if $n_i=n_{i-1}$ 
(resp., $n_i\neq n_{i-1}$).  The resulting skew Young diagram $\ka$
is a border strip and, clearly, 
\eqn\eqCBah{
|\ka| ~=~ M \,,\qquad\qquad \sum_i (s-i) b_{i} ~=~ \sum_i n_i\,,
}
provided $\# \{ i | n_i =k\}\neq 0$ for all $0\leq k\leq n_M$.
This suggests that to each such sequence is associated a $Y(\sln)$
highest weight vector whose expression in terms of spinon states
has a leading term (in general there will be subleading terms with
respect to some ordering on $\{n_i\}$) given by \eqCBaf.  This is
indeed the case [\BS].
The $\sln$ content of the corresponding irreducible $Y(\sln)$ 
module is a subquotient of $L(\La_1)^{\otimes M}$ which can be
understood in terms of a generalized exclusion principle for the spinon
$\ph^i(z)$.  Note that to satisfy the condition on the modes $\{n_i\}$,
one has to insert singlet combinations of $n$ spinons into the `gaps'
between the $n_i$'s.  Hence, the spinon number $M$ in the description 
of \eqCBah\ is not the same as the spinon number $N$ in the corresponding 
result \eqCBad\ for $\sltw$.  See [\BS] for more details, in particular 
on the correspondence between a mode sequence $\{n_i\}$ in \eqCBaf\ and 
a motif, or equivalently a semi-infinite border strip (see (A.17)), 
characterizing an ireducible $Y(\sln)$ highest weight module.


\appendix{A}{Partitions and symmetric functions}

In this appendix we recall the definition of the (skew) Schur functions
associated to a partition (see [\MD] for additional background).

A partition $\la=(\la_1,\la_2,\ldots,\la_r,\ldots)$ is a sequence of
non-increasing, non-negative integers $\la_1\geq \la_2\geq\ldots
\geq 0$ and containing only finitely many non-zero terms.  The length
$l(\la)$ of $\la$ is the number of non-zero elements in $\la$, and the 
weight $|\la|$ of $\la$ is defined by $|\la| = \sum_i \la_i$.
If $N=|\la|$ we say that $\la$ is a partition of $N$.  The set of all
partitions of $N$ is denoted by $\cP_N$, and the set of all
partitions by $\cP$.
For a partition $\la\in\cP_N$, let 
\eqn\eqApaaa{
m_i^\la ~=~ \#\{ j: \la_j=i\} 
}
denote the number of parts of $\la$ that are equal to $i$.  We have
\eqn\eqApaab{
|\la| ~=~ \sum_{i\geq1} \ i m_i^\la\,,\qquad 
l(\la) ~=~ \sum_{i\geq1} \  m_i^\la\,.
}
Sometimes we also denote
\eqn\eqApaac{
\la ~=~ (1^{m_1} 2^{m_2} \ldots r^{m_r} \ldots)\,,
}
where $m_i \equiv m_i^\la$.

As usual, to each partition $\la$ we associate a Young diagram, which we 
also denote by $\la$.  The partition $\la'$, conjugate to $\la$, is
then defined by transposing the Young diagram of $\la$ along the 
main diagonal.  Note that 
\eqn\eqApaad{
\la_i' ~=~ \#\{ j: \la_j\geq i\} \,,
}
such that, in particular, $l(\la) = \la_1'$.  Moreover
\eqn\eqApaae{
m_i^\la ~=~ \la_i' - \la'_{i+1} \,.
}

For any pair of partitions $\la$ and $\mu$, we write $\la\supset\mu$
if $\la_i\geq \mu_i$ for all $i$.  If $\la\supset\mu$ then the 
Young diagram of $\la$ contains the Young diagram of $\mu$.  The
skew Young diagram $\la/\mu$ is then obtained as the set-theoretic 
difference $\la-\mu$.  We say that ${\rm rank}(\la/\mu)=n$ if the 
length of any column of $\la/\mu$ does not exceed $n$, and we put
$|\la/\mu| = |\la| - |\mu|$.

A path in a skew diagram $\la/\mu$ is a sequence $x_0,x_1,\ldots,x_r$
of squares such that $x_{i-1}$ and $x_i$ have a common side, i.e.\
are adjacent.  The diagram $\la/\mu$ is said to be connected if 
any two squares in $\la/\mu$ can be connected by a path.  Finally, a
skew diagram $\la/\mu$ is called a border strip if $\la/\mu$ is connected 
and contains no $2\times2$ block of squares.  We denote the set of
border strips of rank $n$ by ${\rm BS}_n$.

For given $n\in\NN$, 
a semi-standard tableau of shape $\la/\mu$ is an inscription of
the numbers $1,2,\ldots,n$ in each of the boxes of a given skew
Young diagram $\la/\mu$ such that if $a$ and $b$ are the inscriptions
in any pair of adjacent boxes, then
\item{i.} $a<b$ if $b$ is lower-adjacent to $a$,
\item{ii.} $a\geq b$ if $b$ is left-adjacent to $a$.

We denote the set of semi-standard tableaux of shape $\la/\mu$ by
${\rm SST}(\la/\mu)$.  For ${T}\in {\rm SST}(\la/\mu)$ we define 
$m_a(T) = \#\{ a: a\in {T}\}$.

For a skew diagram $\la/\mu$ the skew Schur function $s_{\la/\mu}(x)$
is now defined as 
\eqn\eqApaba{
s_{\la/\mu}(x) ~\equiv~ s_{\la/\mu}(x_1,x_2,\ldots,x_n)
  ~=~ \sum_{{T}\in {\rm SST}(\la/\mu)} \left( \prod_i x_i^{m_i(T)}\right)
  \,.
}
Note that the (normal) Schur functions are included in this definition 
by taking $\mu = \emptyset$.

We have the following expressions for skew Schur functions in terms
of the elementary symmetric functions $e_m(x)$ and $h_m(x)$ [\MD]
\eqn\eqApabb{ \eqalign{
s_{\la/\mu}(x) & ~=~ {\rm det}(h_{\la_i - \mu_j -i+j}(x))_{1\leq i,j\leq r}
  \,,\qquad ({\rm for\ } r\geq l(\la))\cr
& ~=~ {\rm det}(e_{\la'_i - \mu'_j -i+j}(x))_{1\leq i,j\leq s}
  \,,\qquad ({\rm for\ } s\geq l(\la'))  \,.\cr}
}
Note that, in particular,
\eqn\eqApabc{
h_m(x) ~=~ s_{(m)}(x)\,,\qquad e_m(x) ~=~ s_{(1^m)}(x)\,.
}
The skew Schur functions can be expressed in terms of standard Schur 
functions by means of the Littlewood-Richardson rule
\eqn\eqApaxa{
s_{\la/\mu}(x) ~=~ \sum_\nu c^\nu_{\la\mu} s_\nu(x)\,.
}

Now, let $\la/\mu \in {\rm BS}_n$. Let $r$ and $s$ be the number of rows
and columns  of $\la/\mu$, respectively.  Denote the length of the $i$-th
row and column by $a_i$ and $b_i$, respectively.  Then
\eqn\eqApabd{ \eqalign{
a_i & ~=~ \la_i - \mu_i\,, \cr
b_i & ~=~ \la'_{s-i+1} - \mu'_{s-i+1} \,.\cr}
}
We denote $\la/\mu = \langle a_1,\ldots,a_r\rangle = [b_1,\ldots,b_s]$.
Then, \eqApabb\ gives
\eqn\eqApabe{ \eqalign{
s_{\la/\mu}(x) ~=~ s_{\langle a_1,\ldots,a_r\rangle}(x) 
& ~=~ \left| \matrix{
  h_{a_1} & h_{a_1+a_2} && \ldots && h_{a_1+\ldots+a_r} \cr
     1    & h_{a_2}     &&&&\cr
     0    &   1         &&&&\cr
          &   0         &&&& \vdots\cr
          &             & \ddots & \ddots && \cr
     && 0 & 1 & h_{a_{r-1}} & h_{a_{r-1} + a_r} \cr
     &&&  0 & 1 & h_{a_r} \cr} \right| \cr
& \cr
~=~ s_{[b_1,\ldots,b_s]}(x) & ~=~
\left| \matrix{
  e_{b_s} & e_{b_{s} + b_{s-1}} && \ldots && e_{b_s+\ldots+ b_1} \cr
     1    & e_{b_{s-1}}     &&&&\cr
     0    &   1         &&&&\cr
          &   0         &&&& \vdots\cr
          &             & \ddots & \ddots && \cr
     && 0 & 1 & e_{b_2} & e_{b_2+b_1} \cr
     &&&  0 & 1 & e_{b_1} \cr} \right| \,.\cr}
}
Using that
\eqn\eqApabea{
h_a(x) ~=~ \sum_{k_1,\ldots,k_n\geq0\atop \sum k_i =a} \, 
 x_1^{k_1} \ldots x_n^{k_n} \,,
}
we find, for $n=2$ and $a_1,a_2\geq1$,
\eqn\eqApabeb{
h_{a_1}h_{a_2} - h_{a_1+a_2} ~=~ h_{a_1-1}h_{a_2-1}\,.
}
This gives the following drastic simplification of \eqApabe\ for 
$n=2$ (and $a_1,a_r\geq 1$, $a_2,\ldots,a_{r-1}\geq2$)
\eqn\eqApabec{
s_{\langle a_1,\ldots,a_r\rangle} ~=~ h_{a_1-1}h_{a_2-2}\ldots
  h_{a_{r-1}-2} h_{a_r-1} \,.
}
Note furthermore that, e.g.\ because of \eqApaba, we have 
\eqn\eqApabf{
s_{[b_1,\ldots,b_s,n]}(x) ~=~ s_{[b_1,\ldots,b_s]}(x)\,.
}
We may thus extend the definition of skew Schur functions to 
semi-infinite border strips, i.e.\ border strips stabilizing on
columns of length $n$, by
\eqn\eqApabg{
s_{[b_1,\ldots,b_s,n,n,n,\ldots]}(x) ~\equiv~ s_{[b_1,\ldots,b_s]}(x)\,.
}

It is convenient to introduce yet two other parametrizations of ${\rm BS}_n$
which have their origin in integrable spin chain models.
The first parametrization is through a semi-infinite sequence of rapidities
$\{ m_s\}$, i.e.\ a sequence of distinct 
positive integers $0<m_1<m_2<m_3<\ldots$
such that there are no more than $n-1$ consecutive $m_i$'s and such that 
the sequence stabilizes on the sequence $\{\ldots, a+1,a+2,\ldots,
a+(n-1),a+(n+1),\ldots,a+(2n-1),a+(2n+1),\ldots\}$ for some $a$.
Denote the set of such semi-infinite rapidity sequences by $\cR_n^{\infty/2}$.
If $a\equiv k\,{\rm mod\,}n$ for some $0\leq k\leq n-1$ then we have a 
rapidity sequence of conjugacy class $k$.  The vacuum rapidity sequence 
$\{ m_s^{(k)}\}$, in conjugacy class $k$, is given by
$\{ m_s^{(k)}\} = \{ 1,2,3,\ldots,k-1,k+1,\ldots,k+n-1,k+n+1,\ldots\}$.
The second parametrization is through a semi-infinite motif of rank $n$, i.e.\
a sequence $(d_1,d_2,d_3,\ldots)$ with $d_i\in\{0,1\}$, such that there
are no more than $n-1$ consecutive $1$'s and such that 
the sequence stabilizes on the sequence 
$(\ldots,(0,\underbrace{1,\ldots,1}_{n-1})^\infty)$.  We denote the set of 
such motifs by $\cM_n^{\infty/2}$.  
Clearly, we have an isomorphism $\cR_n^{\infty/2}
\simeq
\cM_n^{\infty/2}$ by letting $d_i=1$ iff $i\in\{m_s\}$ and $d_i=0$ otherwise.
We also have $\cR_n^{\infty/2} \simeq {\rm BS}_n$ by associating a 
semi-infinite
border strip
$\langle a_1,a_2,\ldots\rangle \in {\rm BS}_n$, defined by
\eqn\eqApabh{
a_i ~=~ m_i - m_{i-1} \,,\qquad i=1,2,\ldots\,,
}
where $m_0\equiv 0$ with $\{m_s\}\in\cR_n^{\infty/2}$.

Alternatively, given a motif $(d_1,d_2,\ldots) \in \cM_n^{\infty/2}$, 
the corresponding 
border strip $\ka\in{\rm BS}_n$ is constructed as follows:
Write down a square, call this the $0$-th square.  Now, by induction,
write the $i$-th square ($i\in\NN$) under (resp., left to) the
$i-1$-th square if $d_i=1$ (resp., $d_i=0$).

\appendix{B}{The Lie algebras $\sln$ and $\hsln$}

In this appendix we collect some results regarding the Lie algebra 
$\sln$ and its (untwisted) affinization $\hsln$.  We will be rather
brief, our main purpose is to explain the notation used throughout the 
paper.  

The Lie algebra $\gln$ is the algebra of $n\times n$ (complex) matrices.
A basis of $\gln$ is given by the matrix units $e_{ij}, \, 1\leq i,j\leq n$,
with components $(e_{ij})_{kl} = \de_{ik} \de_{jl}$ and commutators
\eqn\eqApba{
[e_{ij} , e_{kl} ] ~=~ \de_{jk} e_{il} - \de_{li} e_{kj} \,.
}
The Cartan subalgebra $\bfh$ of $\gln$ is 
spanned by the $e_{ii},\, 1\leq i\leq n$,
i.e.\ $h\in\bfh$ if and only if $h = \sum a_i e_{ii}$.  Let $\tilde\ep_i,\,
i=1,\ldots,n$, be the orthonormal basis of $\bfh^*$ defined by 
$\tilde\ep_i(e_{jj}) = \de_{ij}$.

Now, let $V$ be any $\gln$ module.
Since the Cartan subalgebra is semi-simple, it acts diagonally
on $V$.  We thus have a decomposition of $V$ as a direct 
sum 
\eqn\eqApbda{
V ~\simeq~ \bigoplus_{\la\in\bfh^*}\, V_\la\,,
}
of eigenspaces $V_\la = \{ v\in V | h\cdot v = \la(h) v\}$.  We call
$\la$ a weight if $V_\la \neq \emptyset$ and we call $V_\la$ a
weight space.

The (abstract) character ${\rm ch}_V$ of a finite-dimensional 
$\gln$ module $V$ is defined as
\eqn\eqApbba{
{\rm ch}_V ~=~ \sum\ {\rm dim}(V_\la)\, e^\la\,.
}
A vector $v\in V$ is called a singular vector 
of weight $\la$ when $v\in V_\la$ and $e_{ij}\cdot v =0$ for all
$i<j$.  A highest weight module $V$ is a module which possesses a
singular vector $v$ such that $V = U(\gln)\cdot v$.  All 
irreducible highest weight modules are therefore 
uniquely characterized by the weight $\la$ of the highest weight vector.
We will denote the irreducible highest weight module with highest weight 
$\la$ by $L(\la)$.
We can expand $\la = \sum_i \la_i \tilde\ep_i$, i.e.\ $\la_i = 
(\la,\tilde\ep_i)$.  The module $L(\la)$ is finite-dimensional if and 
only if $\la_i - \la_{i+1}\in\ZZ_{\geq 0}$ for all $i=1,\ldots,n-1$.

The Weyl group $W(\gln) \simeq S_n$ of $\gln$ acts by permutations
on the $\tilde\ep_i$ thus, upon defining $x_i=e^{\tilde\ep_i}$,
we find that the character ${\rm ch}_V(x)$ is a symmetric function in 
the $x_i$.

The Lie algebra $\sln$ is the subalgebra of $\gln$ consisting of 
traceless $n\times n$ matrices.  
A Chevalley basis for $\sln$, i.e.\ a set of generators satisfying 
\eqn\eqbca{ \eqalign{
[h_i,h_j] & ~=~ 0 \,,\cr
[h_i,x_j^\pm] & ~=~ \pm a_{ij} x_j^\pm \,,\cr
[x_i^+ ,x_j^-] & ~=~ \de_{ij} h_i \,,\cr
({\rm ad}\, x_i^\pm)^{1-a_{ij}} \, x_j^\pm & ~=~ 0\,,\cr}
}
is given by 
\eqn\eqbcb{
x_i^+ ~=~ e_{i i+1} \,,\qquad x_i^- ~=~ e_{i+1 i}\,,\qquad
h_i~=~ e_{ii} - e_{i+1 i+1}\,,\qquad i=1,\ldots,n-1\,.
}
Here $a_{ij} = 2\de_{ij}- \de_{i-1 j} - \de_{i+1 j}$ is the Cartan 
matrix of $\sln$. 

An (overcomplete) basis for the dual Cartan subalgebra $\bfh^*$ of $\sln$ is 
given by $\ep_i,\,i=1,\ldots,n$, defined by $\ep_i(\sum_j a_j e_{jj}) = a_i$.
They satisfy the constraint $\sum_i \ep_i =0$.  Note that in terms of the
corresponding basis for $\gln$ we have
\eqn\eqApbca{
\ep_i ~=~ \tilde\ep_i - \textstyle{1\over n} \,\sum_j \,\tilde\ep_j\,.
}
Consequently,
\eqn\eqApbcb{
(\ep_i,\ep_j) ~=~ \de_{ij} - \textstyle{1\over n} \,.
}
In terms of the $\ep_i$, the simple roots $\al_i$ of $\sln$ are given by 
\eqn\eqApbcc{
\al_i ~=~ \ep_i - \ep_{i+1}\,,\qquad i=1,\ldots,n-1\,,
}
while the fundamental weights $\La_i$, defined by $(\La_i,\al_j)=\de_{ij}$,
are given by 
\eqn\eqApbcd{
\La_k ~=~ \ep_1 + \ep_2 + \ldots + \ep_k\,,\qquad k=1,\ldots,n-1\,.
}
Furthermore, the inverse Cartan matrix is given by
\eqn\eqApbce{
(a^{-1})_{ij} ~=~ (\La_i,\La_j) ~=~ {\rm min}(i,j) - \textstyle{1\over n} \,
 (ij)\,.
}
In particular
\eqn\eqApbcea{
\half |\La_k|^2 ~=~ { k(n-k)\over 2n }\,.
}

The weight lattice, i.e.\ the lattice spanned by the fundamental weights
$\La_i$ is denoted by $P$, while the set of weights in $P$ of
conjugacy class $k$, i.e.\ those weights $\sum m_i\La_i$ such that
$\sum im_i \equiv k\ {\rm mod}\ n$ is denoted by $P^{(k)}$.  The
root lattice $P^{(0)}$ is also denoted by $Q$.

The irreducible highest weight modules $L(\la)$ are labeled by
a highest weight vector $\la = \sum_i m_i \La_i$ and are
finite-dimensional if and only if $m_i\in\ZZ_{\geq0}$, i.e.\ iff
$\la$ is a dominant integral weight.  The set of dominant integral
weights is denoted by $P_+$.  The 
irreducible finite-dimensional modules are thus in 1--1
correspondence with Young diagrams of rank $n-1$, i.e.\ 
$\la = (\la_1,\la_2,\ldots,\la_{n-1})$ where 
$\la_i = m_i + m_{i+1} + \ldots + m_{n-1}$.
Equivalently, $\la_i = (\la,\ep_i-\ep_n)$.
Also note that the $\ep_i,\,i=1,\ldots,n$ are the weights of
the $n$-dimensional fundamental irreducible representation $L(\La_1)$
of $\sln$.

Again, upon identifying $x_i = e^{\ep_i}$, the character of $L(\la)$
becomes a symmetric function in the $x_i$ (here it is implicitly 
understood that $x_1\ldots x_n=1$ because of the constraint 
$\sum \ep_i=0$).  In fact,
\eqn\eqApbcf{
{\rm ch}_{L(\la)}(x) ~=~ s_\la(x)\,.
}
In particular we have 
\eqn\eqApbcg{
{\rm ch}_{L(\La_m)}(x) ~=~ e_m(x) \,,\qquad
{\rm ch}_{L(m\La_1)}(x)  ~=~ h_m(x) \,.
}

For any simple Lie algebra $\bfg$ (in this paper we will only
consider $\bfg \simeq \sln$), the untwisted affine Lie algebra 
$\whg$ is defined as a central extension of the loop algebra 
$\bfg \otimes \CC[t^{-1},t]$, i.e.\
\eqn\eqApbea{
\whg ~=~ ( \bfg \otimes \CC[t^{-1},t] ) \oplus \CC \,\ell \oplus L_0\,,
}
where $\ell$ is the central element and $L_0$ a derivation.  To be 
precise, for $x,y\in\bfg$, $m,n\in\ZZ$,
\eqn\eqApbeb{ \eqalign{
[x\otimes t^m , y\otimes t^n ] & ~=~ [x,y] \otimes t^{m+n} + m \,\ell\,
\de_{m+n,0}
  (x,y) \,,\cr
[L_0, x\otimes t^m] & ~=~ -m (x\otimes t^m) \,,\cr
[\ell, x\otimes t^m] & ~=~ [\ell, L_0] ~=~ 0 \,.\cr}
}
Henceforth we will write $x_m$ for the element $x\otimes t^m$ of $\whg$.
The subalgebra consisting of the elements $x_0=x\otimes t^0$ is identified
with $\bfg$.  In the remainder we consider $\bfg=\sln$.
The fundamental weights of $\hsln$ are denoted by $\widehat{\La}_i$,
$i=0,\ldots,n-1$

In each irreducible $\whg$ module $V$, the central element $\ell$ acts by
a constant (also to be denoted by $\ell$) which is referred to as the 
level of $V$.  

In this review we consider two important classes of modules of 
$(\hsln)_\ell$.  For $\ell=0$ we have the so-called evaluation modules
$L(\La)(\!(z)\!) \equiv {\rm ev}_z^* (L(\La))$ defined by pulling
back a finite dimensional  irreducible $\sln$ module $L(\La)$ by means
of the evaluation homomorphism ${\rm ev}_z : U(\sln \otimes \CC[t^{-1},t])
\to U(\sln)$ defined by ${\rm ev}_z(x_m) = z^m x$.

The other important class of modules are the so-called 
integrable highest weight modules.  They exists for $\ell\in\NN$ and 
highest weights $\widehat{\La} = \sum_{i=0}^{n-1} m_i \widehat{\La}_i$,
$m_i\in\ZZ_{\geq0}$,  such that $\sum m_i =\ell$.  Let us denote the
set of such weights by $\widehat{P}^{\ell}_+$. 

The weight space decomposition of an integrable highest weight 
module $V$ with highest weight $\widehat{\La}\in \widehat{P}^{\ell}_+$ 
is now defined 
with respect to $\bfh$, the Cartan subalgebra of $\sln$, and $L_0$.
That is, $V\simeq \oplus V_{\la,n}$ where 
\eqn\eqApbed{
V_{\la,n} ~=~ \{ v\in V\,|\, h\cdot v = \la(h)\,v,\ L_0\cdot v = 
  (\De(\La) + n) v\}\,.
}
where 
\eqn\eqApbee{
\De(\widehat{\La}) ~=~ { (\La,\La+2\rh) \over 2(\ell+n) }\,
}
is the conformal dimension of the highest weight vector of 
$L(\widehat{\La})$.
Here, $\La$ is the projection of the $\hsln$ weight $\widehat{\La}$
onto $\bfh^*$ and $\rh$ is the Weyl vector of $\sln$.
In particular, for $\ell=1$, we have (cf.\ \eqApbcea)
\eqn\eqApbeea{
\De_k ~\equiv~ \De(\widehat{\La}_k) ~=~ { k(n-k)\over2n}\,.
}

The character of the integrable highest weight module 
$L(\widehat{\La})$ is defined as
\eqn\eqApbf{
{\rm ch}_{L(\widehat{\La})} ~=~ q^{\De(\widehat{\La})} \sum_{\la,n} 
{\rm dim}(L(\widehat{\La})_{\la,n}) \, q^n\, e^\la\,.
}

The characters of the level-$1$ integrable highest weight modules 
$L(\widehat{\La}_k),\, k=0,\ldots,n-1$ 
of $\hsln$ are given by
\eqn\eqApbec{
{\rm ch}_{L(\widehat{\La}_k)} ~=~ \sum_{\la\in P^{(k)}} {q^{{1\over2}|\la|^2}
  \over (q)_\infty^{n-1} } e^\la\,,
}
or, more explicitly,
\eqn\eqApbeca{
{\rm ch}_{L(\widehat{\La}_k)}(x;q) ~=~ \sum_{ k_1,\ldots,k_n \in
  \ZZ - {k\over n} \atop \sum k_i =0} {q^{ {1\over2} (k_1^2 +\ldots+
  k_n^2)} \over (q)_\infty^{n-1}} x_1^{k_1} \ldots x_n^{k_n}\,.
}

Chiral vertex operators (CVO's) at level-$\ell$ are $\hsln$ intertwiners
\eqn\eqApbebg{
\Ph \left( \matrix{ \la_3 \cr \la_2 \ \la_1 \cr} \right) \ 
  : \ L(\widehat{\la}_1) \otimes L(\la_3)(\!(z)\!) ~\to~
  L(\widehat{\la}_2)\,,
}
where $\la_3\in P_+$ and $\widehat{\la}_i\in \widehat{P}^\ell_+$ for $i=1,2$.  
It is convenient to think of CVO's as a collection of linear maps 
\eqn\eqApbebh{
\Ph^i \left( \matrix{ \la_3 \cr \la_2 \ \la_1 \cr} \right) (z) \ : \ 
  L(\widehat{\la}_1)~\to~ L(\widehat{\la}_2)\,,
}
transforming under $(\hsln)_0$ as
\eqn\eqApbebi{
[x_m , \Ph^i \left( \matrix{ \la_3 \cr \la_2 \ \la_1 \cr} \right) (z)]
  ~=~ z^m  \rh(x)^i{}_j\cdot 
  \Ph^j \left( \matrix{ \la_3 \cr \la_2 \ \la_1 \cr} \right) (z)\,,
}
where $\rh(x)^i{}_j$ is 
the action of $x$ in the $\sln$ irreducible representation
$L(\la_3)$.
A CVO $\Ph \left( \matrix{ \la_3 \cr \la_2 \ \la_1 \cr} \right)$ exists
if and only if $\la_2$ occurs in the fusion rule $\la_1 \times \la_3$.
The CVO's for $\la_3 = \La_1$ will be referred to a spinon operators.  
Since the level-$1$ fusion rules read $\La_k \otimes \La_1 = \La_{k+1}$
(here the subscripts are taken modulo $n$) the only non-vanishing
spinon operators at $\ell=1$ are
\eqn\eqApbebj{
\Ph^i \left( \matrix{ \La_1 \cr \La_{k+1} \ \La_k \cr} \right) (z)\,.
}


\appendix{C}{Yangians and their representations}

Yangians were introduced by Drinfel'd [\Dra,\Drb] (see [\Drc,\CPc] for
reviews).  Here we will briefly review their definition and some
aspects of their representation theory.  For definiteness we will 
restrict ourselves to the algebras $\sln$ and $\gln$, the cases 
of interest to this paper, but most of what follows applies to
arbitrary $\bfg$ as well.

The Yangian $Y(\bfg)$ is defined as the (unique) flat deformation
of the loop algebra $U(\bfg[t])$ in the class of Hopf algebras.
There are three equivalent realizations, each useful for 
different purposes.
\thm\thApca
\proclaim Definition \thApca\ [\Dra].  The Yangian $Y(\gln)$ is defined to be 
the Hopf algebra generated by elements $t^{(p)}_{ij}$,
$i,j=1,\ldots,n,\, p\in\ZZ_{\geq0}$, with relations
\eqn\eqApca{
[t^{(p+1)}_{ij}, t^{(q)}_{kl}] - [t^{(p)}_{ij}, t^{(q+1)}_{kl}] ~=~
  - ( t^{(p)}_{kj} t^{(q)}_{il} - t^{(q)}_{kj} t^{(p)}_{il} )\,,
}
where we have defined $t^{(-1)}_{ij} = \de_{ij}$.\par

Let $t_{ij}(u) \in Y(\gln)[[ u^{-1}]]$ be defined by $t_{ij}(u)
= \de_{ij} + \sum_{p\geq0} t^{(p)}_{ij} u^{-p-1}$.  Then,
after introducing an `auxiliary space' $\CC^n$, and writing 
$t(u) = \sum_{i,j} t_{ij}(u) \otimes e_{ij}$, where $e_{ij}$ is 
the matrix with
components $(e_{ij})_{kl} = \de_{ik}\de_{jl}$, we can also write
\eqApca\ as
\eqn\eqApcb{
R(u-v) (t(u) \otimes 1)(1\otimes t(v)) ~=~ (1\otimes t(v))
  (t(u) \otimes 1)R(u-v) \,,
} 
where $R(u) = 1 + {1\over u} \sum_{i,j} e_{ij} \otimes e_{ji} \in
{\rm End}(\CC^n\otimes\CC^n)$.  
The matrix $R(u)$ corresponds to the simplest (rational) solution of
the Yang-Baxter equation.

The co-multiplication $\De$, co-unit $\ep$, and anti-pode $S$ are given by
\eqn\eqApcc{ \eqalign{
\De(t_{ij}(u)) & ~=~ \sum_k  t_{ik}(u) \otimes t_{kj}(u)\,,\cr
\ep(t_{ij}(u)) & ~=~ \de_{ij}\,,\cr
S(t(u)) & ~=~ t(u)^{-1} \,.\cr}
}
\thm\thApcb
\proclaim Theorem \thApcb\ [\Dra].  
The center $Z(Y(\gln))$ of $Y(\gln)$ is generated 
by the quantum determinant 
\eqn\eqApcd{ 
{\rm det}_q\, t(u) 
~=~ \sum_{\pi \in S_n}
  {\rm sgn}(\pi) t_{1\pi(1)}(u + {n-1\over2})
  \ldots t_{ n\pi(n)}(u - {n-1\over2}) 
}
and 
\eqn\eqApce{
\De( {\rm det}_q\, t(u) ) ~=~ {\rm det}_q\, t(u) \otimes {\rm det}_q\, t(u)\,.
}
\par

The theorem shows that we can define $Y(\sln)$ as the Hopf subalgebra
of $Y(\gln)$ consisting out of those $t(u)$ with ${\rm det}_q\, t(u)=1$.

Note that we have a family of automorphisms  $\ta_a\,:\, Y(\gln)
\to Y(\gln)$ ($a\in\CC$) defined by 
\eqn\eqApcf{
\ta_a (t_{ij}(u)) ~=~ t_{ij}(u-a)\,.
}

The second realization of $Y(\sln)$ is given by
\thm\thApcc
\proclaim Definition \thApcc\ [\Drb].  The Yangian $Y(\sln)$ is defined to be
the Hopf algebra with generators 
$\{ x^{\pm}_{ik}, h_{ik}\}$, $i=1,\ldots,n-1,\, k\in\ZZ_{\geq0}$,
and relations
\eqn\eqApcg{
[h_{ik},h_{jl}]~=~ 0\,,\qquad 
[h_{i0},x^{\pm}_{jl}] ~=~ \pm a_{ij} x^{\pm}_{jl} \,,\qquad
[x^+_{ik} , x^-_{jl}] ~=~ \de_{ij} h_{i k+l}\,,
}
\eqn\eqApch{
[ h_{i k+1}, x^{\pm}_{jl} ] - [h_{ik},x^{\pm}_{j l+1} ] ~=~
  \pm{\half} a_{ij} ( h_{ik} x^{\pm}_{jl} + x^{\pm}_{jl}h_{ik} ) \,,
}
\eqn\eqApci{
[ x^{\pm}_{i k+1}, x^{\pm}_{jl} ] - [x^{\pm}_{ik},x^{\pm}_{j l+1} ] ~=~
  \pm{\half} a_{ij} ( x^{\pm}_{ik}x^{\pm}_{jl}+ x^{\pm}_{jl}x^{\pm}_{ik})\,,
}
\eqn\eqApcj{
\sum_{\pi\in S_m} [x^{\pm}_{i k_{\pi(1)}}, [x^{\pm}_{i k_{\pi(2)}}, \ldots
  [x^{\pm}_{i k_{\pi(m)}}, x^\pm_{j l}] \ldots ]] ~=~ 0 \,,\qquad
  {\rm where\ } m=1-a_{ij},\,i\neq j\,.
}
where $a_{ij}$ is the Cartan matrix of $\sln$.\par

Again, it is sometimes convenient to work with the generating series
\eqn\eqApcm{
h_i(u) ~=~ 1 + \sum_{k\geq0} h_{ik} u^{-k-1}\,,\qquad
x^\pm_i(u) ~=~ \sum_{k\geq0} x^\pm_{ik} u^{-k-1}\,.
}
For example, a family of Hopf algebra automorpisms $\ta_a,\, a\in\CC$,
is defined by 
\eqn\eqApcma{
\ta_a(x_i^\pm(u)) ~=~ x_i^\pm(u-a)\,,\qquad 
\ta_a(h_i(u)) ~=~ h_i(u-a)\,.
}

In order to exhibit the isomorphism between the realizations \thApca\ and 
\thApcc, define 
\eqn\eqApck{ \eqalign{
a_i(u) & ~=~ ({\rm det}_q\,t_{rs}(u))_{r,s\in\{1,\ldots,i\}} \,, \qquad 
  a_0(u) ~=~ 1\,,\cr
b_i(u) & ~=~ ({\rm det}_q\,t_{rs}(u))_{r\in\{1,\ldots,i\},s\in
  \{1,\ldots, i-1,i+1\}} \,,\cr
c_i(u) & ~=~ ({\rm det}_q\,t_{rs}(u))_{r\in\{1,\ldots, i-1,i+1\},s\in
  \{1,\ldots,i\}} \,,\cr}
}
where ${\rm det}_q\,t(u)$ is the quantum determinant defined in \eqApcd.
Then we have
\thm\thApce
\proclaim Theorem \thApce\ [\Drb].  A (non-canonical) isomorphism between 
the definitions \thApca\ and \thApcc\ of $Y(\sln)$ is given by
\eqn\eqApcl{ \eqalign{
x_i^+(u) & ~=~ b_i(u) a_i(u)^{-1} \,,\cr
x_i^-(u) & ~=~ a_i(u)^{-1} c_i(u) \,,\cr
h_i(u)   & ~=~ a_i(u)^{-1}a_i(u+1)^{-1} a_{i-1}(u+\half)
               a_{i+1}(u+\half) \,.\cr}
}
\par

Finally, the third realization of $Y(\sln)$ is
\thm\thApcd
\proclaim Definition \thApcd\ [\Dra]. 
The Yangian $Y(\sln)$ is defined to be the 
Hopf algebra generated by $x, J(x)$, $x\in\sln$, where $J$ is a linear
functional on $\sln$, with relations
$$
[x,J(y)] ~=~ J( [x,y]) \,,
$$
$$
[J(x), J([y,z]) ] + {\rm cycl.} ~=~ {\textstyle{1\over4}} 
  \sum_{a,b,c} ([x,I_a], [[y,I_b],[z,I_c]]) \{ I_a,I_b,I_c \} \,,
$$
$$ \eqalign{
& [[J(x),J(y)],[z,J(w)]] + [[J(z),J(w)],[x,J(y)]] \cr
  & ~=~ {\textstyle{1\over4}} \sum_{a,b,c} \big( 
  ([x,I_a],[[y,I_b],[[z,w],I_c]]) + ([z,I_a],[[w,I_b],[[x,y],I_c]])\big)
  \{ I_a,I_b,J(I_c) \} \,.\cr}
$$
where $\{I_a\}$ denotes
an orthonormal basis of $\sln$ with respect to the 
Killing form $(\ ,\ )$, and 
$$
\{ x_1,\ldots, x_n\} ~=~ {1\over n!} \sum_{\pi\in S_n} x_{\pi(1)}
  \ldots x_{\pi(n)}\,,
$$
denotes the symmetrizer.

Let us now discuss the representations of $Y(\sln)$.  We will be interested
in the finite-dimensional representations in particular.  
\thm\thApcaa
\proclaim Definition \thApcaa. A representation $V$ of $Y(\sln)$ is said to
be a highest weight module if it is generated by a vector 
$v\in V$, i.e.\ $V = Y(\sln)\cdot v$, such that 
\eqn\eqApcaa{
x^{+}_{ik} \cdot v ~=~ 0 \,,\qquad h_{ik}\cdot v ~=~ d_{ik} v\,,
}
for all $i=1,\ldots, n-1$, $k\in\ZZ_{\geq0}$ and some $\underline{d}=
d_{ik}\in\CC$.\par

As usual, the Verma module $M(\underline{d})$ is universal in the class
of highest weight modules, i.e.\ every highest weight module is a 
quotient module of the Verma module, and $M(\underline{d})$ has a unique 
irreducible quotient $L(\underline{d})$.

For a highest weight module $V(\underline{d})$ define $d_i(u)\in\CC[[u^{-1}]]$
by $h_i(u)\cdot v = d_i(u) v$, i.e.\ 
\eqn\eqApcab{
d_i(u) ~=~ 1 + \sum_{k\geq0} d_{ik} u^{-k-1}\,.
}
The classification of finite-dimensional representations of $Y(\sln)$ 
is due to Drinfel'd 
\thm\thApcab
\proclaim Theorem \thApcab\ [\Drb].  
\item{(i)} Every finite-dimensional irreducible representation of
$Y(\sln)$ is a highest weight module.
\item{(ii)} The irreducible representation $L(\underline{d})$ is 
finite-dimensional
if and only if there exist polynomials $P_i(u) \in \CC[u]$, $i=1,\ldots,
n-1$, such that 
\eqn\eqApcac{
d_i(u) ~=~ { P_i(u+1) \over P_i(u)} \,.
}\par

\thm\thApcac
\proclaim Corollary \thApcac. There exists a 1--1 correspondence 
between finite dimensional irreducible representations of $Y(\sln)$
and sets of monic polynomials $P_i(u) \in \CC[u]$, $i=1,\ldots,
n-1$.\par

Note that for a finite-dimensional irreducible module $V$ 
associated to Drinfel'd polynomials $P_i(u)$, the module $\ta_a(V)$
is also finite-dimensional, irreducible and the associated polynomials
are $P_i(u-a)$ (see \eqApcma).

Recall that for $U(\bfg[t])$ we can easily construct 
finite-dimensional irreducible representations 
by pulling back
finite-dimensional irreducible representations $V$ of $U(\bfg)$ by means of 
the evaluation homomorphism ${\rm ev}_a : U(\bfg[t]) \to
U(\bfg)$ defined by ${\rm ev}_a(x \otimes t^k) = a^k x$ for some $a\in\CC$.
A representation $V(a) = {\rm ev}_a^* (V)$ of $U(\bfg[t])$ is called 
an evaluation representation.  It turns out, moreover, that every
finite-dimensional irreducible representation of $U(\bfg[t])$ is
isomorphic (as a $U(\bfg)$-module) to a tensor product of evaluation 
representations [\CPb].  Since $Y(\bfg)$ is a deformation of $U(\bfg[t])$
it is natural to ask whether a similar construction of all finite
dimensional irreducibles exists for $Y(\bfg)$.  For general $\bfg$,
an evaluation homomorphism does not exist and, consequently, 
finite-dimensional irreducible representations $V$ 
of $U(\bfg)$ in general do not
extend to representations of $Y(\bfg)$.  It turns out that usually
one can extend $V$ by `smaller' representations such that one can define
an irreducible action of $Y(\bfg)$ on the extension. It is an important 
open problem to find the minimal such extension.  
For $\bfg=\sln$ we have, however,
\thm\thApcad
\proclaim Theorem \thApcad.  Any finite-dimensional representation
$(V,\rh)$ of $\gln$ extends to a representation of $Y(\gln)$. Explicitly,
\eqn\eqApcad{
\rh( t_{ij}(u) ) ~=~ 1 - {1\over u}\, \rh(e_{ij})\,.
}
\par

On a highest weight vector $v\in V$ of an irreducible $Y(\gln)$
highest weight module, $t(u)$ takes the form
\eqn\eqApcae{
t(u)\cdot v ~=~ \left( \matrix{
 t_{11}(u) & 0 & 0 & \ldots & 0 \cr
   * & t_{22}(u) & 0 && 0 \cr
   \vdots &&&& \vdots \cr
  * & *  & \ldots && t_{nn}(u) \cr} \right) \cdot v\,.
}
Thus, it follows from \eqApcl\ that
\eqn\eqApcaf{
d_i(u)  ~=~ { t_{i+1 i+1}(u - {i-1\over2}) \over 
  t_{ii}(u - {i-1\over2})} \,,
}
which provides the connection of the Drinfel'd polynomials \eqApcac\
with the dialgonal eigenvalues of $t(u)$ on the highest weight vector.
In particular, for the irreducible representations $L(\la)$ 
($\la_i= (\la,\tilde\ep_i)$), 
of $Y(\gln)$ given by Theorem \thApcad, we find
\eqn\eqApcag{
d_i(u) ~=~ { u - {i-1\over2} - \la_{i+1} \over u - {i-1\over2} - \la_{i} }
  ~=~ 1 + (\la_{i}- \la_{i+1}) {1\over u} + \cO({1\over u^2})\,,
}
such that 
\eqn\eqApcah{ \eqalign{
P_i(u) & ~=~ \prod_{j=0}^{\la_{i}- \la_{i+1}-1} \, \left(
  u - {i-1\over2} - \la_{i} + j \right)\cr
& ~=~ \prod_{j=1\atop \la_j'=i}^{\la_1} \, \left(
  u - \half\la_j' -j + \half \right)\,.\cr}
}

Note that the roots of the polynomials $P_i(u)$ for the evaluation
representations of Theorem \thApcad\ form a string on the real
axis of length $\la_i - \la_{i+1}$, i.e.\ a set of the form
$\{ a,a+1,\ldots,a+(n-1)\}$ where $n=\la_i - \la_{i+1}$.

For $Y(\sltw)$ we have the following analogue of the abovementioned theorem
for $U(\sltw[t])$
\thm\thApcae
\proclaim Theorem \thApcae\ [\CPa]. Every finite-dimensional representation
of $Y(\sltw)$ is isomorphic to a tensor product of evaluation representations
$L(\la)(a)$.\par

In fact, for a given Drinfel'd polynomial $P(u)$, the corresponding 
(irreducible) tensor product is easily obtained by collecting the 
roots of $P(u)$ in strings and associating an evaluation representation
with each string (see [\CPa] for a more precise statement).

For $Y(\sln),\, n\geq3$ the situation is not that simple.  In fact,
the character of the irreducible module associated to an arbitrary
set of monic polynomials $P_i(u),\, i=1,\ldots,n-1$ (see Corollary
\thApcac) is not known in general.  
Fortunately, for our purposes it suffices to consider a special
subclass, the so-called tame modules [\NT].
They are defined as follows:  Consider the canonical filtration
\eqn\eqApcba{
Y(\frak{gl}_1) ~\subset~ Y(\frak{gl}_2) ~\subset~ \ldots ~\subset~
Y(\gln)\,.
}
The algebra $A(\gln)$, generated by all the centres of all the 
algebras in the chain \eqApcba, is a maximally commutative subalgebra  
of $Y(\gln)$, the so-called Gel'fand-Zetlin algebra.  
A $Y(\gln)$-module $V$ is called tame, if the
subalgebra $A(\gln)$ acts semi-simply on $V$.  In fact, it turns out that
the action of $A(\gln)$ in every irreducible module is simple.  The
eigenbasis of $A(\gln)$ in such an irreducible module is called a
Gel'fand-Zetlin basis.
\thm\thApcba
\proclaim Theorem \thApcba\ [\Ol].  Let $N\in\ZZ_{\geq0}$ be arbitrary.
\item{i.} We have a homomorphism 
$$
\imath\,:\, Y(\gln) ~\rightarrow~ [U(\frak{gl}_{N+n})]^{ \frak{gl}_N }\,,
$$
where the right hand side denotes the commutant of $\frak{gl}_N$ in 
$U(\frak{gl}_{N+n})$.  This homomorphism
becomes injective for $N\to\infty$.
\item{ii.} The image of the homomorphism $\imath$, together with
the centre $Z(\frak{gl}_N)$, generates the commutant.\par

Now, for any dominant integral weights $\la$ and $\mu$
of the Lie algebras $\frak{gl}_{N+n}$ and $\frak{gl}_N$, respectively,
denote by $L(\la,\mu)$ the subspace of all $\frak{gl}_N$ singular vectors
of weight $\mu$ in the irreducible $\frak{gl}_{N+n}$ module 
$L(\la)$ of highest weight $\la$, i.e.\ the multiplicity of $L(\mu)$ in
$L(\la)$ under the canonical embedding $\frak{gl}_N \subset \frak{gl}_{N+n}$.
Note that, in particular, $L(\la,\mu) \neq \emptyset$ implies
$\la\supset\mu$ as an inclusion of Young diagrams.

The homomorphism of Theorem \thApcba\ equips $L(\la,\mu)$ with the 
structure of a $Y(\gln)$ module.  In fact
\thm\thApcbb
\proclaim Theorem \thApcbb\ [\NT].  
\item{i.} $L(\la,\mu)$ is an irreducible tame $Y(\gln)$ module.
\item{ii.} Every irreducible tame $Y(\gln)$ module splits into
a tensor product of modules of the form  $L(\la,\mu)(a) = 
\ta_a(L(\la,\mu))$.\par

In fact, it is possible to explicitly describe the action 
of $Y(\gln)$ on $L(\la,\mu)$ on a basis that diagonalizes
the GZ-algebra $A(\gln)$.  Such a basis is labeled by the 
GZ-schemes of $L(\la,\mu)$.  A GZ-scheme $\backslash\La/$ of
$L(\la,\mu)$ is an array of (non-negative) integers $\la_{m,i},\,
m=0,\ldots,n,\, i=1,\ldots,N+m$, i.e.\
$$
\backslash\La/ ~=~ \matrix{
  \la_{n,1} & & \la_{n,2} & \ldots && \ldots&\la_{n,N+n} \cr
    & \la_{n-1,1} & \ldots && \ldots & \la_{n-1,N+n-1} \cr
    && \la_{0,1} &\ldots & \la_{0,N} \cr}
$$
such that $\la_{n,i}=\la_i$, $\la_{0,i}=\mu_i$ and satisfying the 
condition $\la_{m,i} \geq \la_{m-1,i} \geq \la_{m,i+1}$ for all
$m$ and $i$.  Note that, because of these conditions, we have a 
sequence of partitions 
\eqn\eqApcbc{
\la^{(0)} ~\subset~ \la^{(1)} ~\subset~ \ldots ~\subset~ \la^{(n)}
}
where $\la^{(m)} = (\la_{m,1},\la_{m,2},\ldots,\la_{m,N+m})$.

The $\sln$ weight of the basis vector labeled by $\backslash\La/$
is given by 
\eqn\eqApcbb{
\sum_{m=1}^n \left( \sum_{i=1}^{N+m} \la_{m,i} - \sum_{i=1}^{N+m-1}
  \la_{m-1,i} \right) \ep_m \,.
}
Thus we easily deduce that there exists a weight-preserving, 1--1
correspondence between the GZ-schemes $\backslash\La/$ of 
$L(\la,\mu)$ and the semi-standard tableaux $T\in {\rm SST}(\la/\mu)$.
This proves
\thm\thApcbc
\proclaim Theorem \thApcbc\ [\KKN].  The $\sln$ character of the irreducible 
$Y(\gln)$ module $L(\la,\mu)$ is given by the skew Schur function 
$s_{\la/\mu}(x)$, i.e.\
\eqn\eqApcbe{
{\rm ch}_{L(\la,\mu)}(x) ~=~ s_{\la/\mu}(x)\,.
}\par

In general the irreducible $Y(\sln)$ module $L(\la,\mu)$ is 
reducible under $\sln \subset Y(\sln)$.  The decomposition is given
by \eqApaxa.  {}From the 
explicit action of the Yangian generators [\NT] one furthermore
concludes
\thm\thApcbd
\proclaim Theorem \thApcbd\ [\NT].  The Drinfel'd polynomials of the 
irreducible $Y(\sln)$ modules $L(\la,\mu)$ are given by
\eqn\eqApcbd{
P_i(u) ~=~ \prod_{j=1\atop \la_j'-\mu_j' =i}^{\la_1} \left(
 u - \half ( \la_j' + \mu_j') - j + \half \right)\,.
}
\par

Note that for $\mu=\emptyset$ we recover \eqApcah.


\appendix{D}{$q$-identities}

In this section we recall some elementary $q$-identities.  The proofs
can be found in [\An].

Recall the definition of the $q$-number
\eqn\eqApda{
(z;q)_N ~=~ \prod_{k=1}^N (1 -z q^{k-1})\,.
}
We will write $(q)_N = (q;q)_N$ for short.  We have the following 
useful expansions 
\thm\thApdb
\proclaim Theorem \thApdb.  Let $N\in\NN$, then
\item{i.}
\eqn\eqApdf{ 
(z;q)_N      ~=~ \sum_{n=0}^N\ (-z)^n q^{ {1\over 2}n(n-1)} 
  \qbin{N}{n} \,,
}
\item{ii.}
\eqn\eqApdg{
(z;q)_N^{-1}  ~=~ \sum_{n\geq0} \ \qbin{N+n-1}{n} z^n \,.
}
\par

Here, the $q$-multinomial is defined by
\eqn\eqApdb{
\qbin{k_1+k_2+ \ldots + k_n}{k_1, k_2, \ldots, k_n} ~=~
 { (q)_{k_1+k_2+ \ldots + k_n} \over (q)_{k_1} (q)_{k_2}\ldots (q)_{k_n}}\,,
}
and we have put $\qbin{N}{n,N-n} = \qbin{N}{n}$ for short.
They are the coefficients in the expansion of the generalized 
Rogers-Sz\"ego polynomial in $n$ variables
\eqn\eqApdd{
H^{(n)}_N(x;q) ~=~ \sum_{k_1 + k_2 + \ldots + k_n = N} \ 
  \qbin{k_1+k_2+ \ldots + k_n}{k_1, k_2, \ldots, k_n} x_1^{k_1} \ldots
  x_n^{k_n} \,.
}
The generating function for the Rogers-Sz\"ego polynomial is given by
\eqn\eqApde{
\sum_{N\geq0}\ H^{(n)}_N(x;q) {t^N \over (q)_N} ~=~
  {1\over (tx_1;q)_{\infty} \ldots (tx_n;q)_{\infty} } \,.
}
Note that for $q=1$ we can interpret $H_N^{(n)}(x;q)$ as the character
of the $\sln$ module $L(\La_1)^{\otimes N}$ through the 
identification $x_i = e^{\ep_i}$ (see App.\ B), while $H_N^{(n)}(x;q^{-1})$
for arbitrary $q$ is proportional to the character of the so-called 
$\sln$ Haldane-Shastry spin chain of length $N$ [\Hik,\BS].

A very useful lemma is the following
\thm\thApda
\proclaim Lemma \thApda\ (Durfee square).  For all $m\in\ZZ$ we have 
\eqn\eqApdc{
\sum_{a,b\geq0 \atop a-b= m} { q^{ab} \over (q)_a (q)_b } ~=~
  {1\over (q)_\infty} \,.
}\par

Finally, to derive the spinon form of the $\hsln$ characters, we need
the following lemma 
\thm\thApdc
\proclaim Lemma \thApdc. For all $M,N \in \ZZ_{\geq0}$, we have
\item{i.} 
\eqn\eqApdl{
\sum_{m\geq 0} \ (-1)^m { q^{ {1\over2} m(m-1)} \over
  (q)_{M-m} (q)_{N-m} (q)_m } ~=~ { q^{MN} \over (q)_M (q)_N } \,.
}
\item{ii.}
\eqn\eqApdm{
\sum_{m\geq 0} \ { q^{(M-m)(N-m)} \over (q)_{M-m} (q)_{N-m} (q)_m }
  ~=~ {1\over (q)_M (q)_N }\,.
}\par

\ni {\it Proof:} \hfil\break
i.\ Multiplying the left hand side by $(q)_N z^M$ and summing over $M$, we
find
$$ \eqalign{
\sum_{M\geq0} & \left( \sum_{m\geq0} \ (-1)^m { q^{ {1\over2} m(m-1)}
  \over (q)_{M-m} } \qbin{N}{m} \right) z^M 
  ~=~ \left( \sum_{M\geq0} \ {z^M\over (q)_M } \right)
  \left( \sum_{m\geq0}\ (-z)^m q^{ {1\over2} m(m-1) }\qbin{N}{m} \right) \cr
& ~=~ { (z;q)_N \over (z;q)_\infty } 
 ~=~ { 1 \over (zq^N;q)_\infty } 
 ~=~ \sum_{M\geq0}\ { q^{MN} z^M \over (q)_M } \cr}
$$
ii.\ Again, multiplying the left hand side by $(q)_N z^M$ and 
summing over $M$, we find
$$\eqalign{
\sum_{M\geq0} & \left( \sum_{m\geq0} { q^{(M-m)(N-m)} \over (q)_{M-m} }
 \qbin{N}{m} \right) z^M 
  ~=~ \sum_{M,m\geq0} { q^{M(N-m)} \over (q)_M} \qbin{N}{m} z^{M+m} \cr
& ~=~ \sum_{m\geq0} {1\over (zq^{N-m};q)_\infty} \qbin{N}{m} z^m 
  ~=~ {1\over (z;q)_\infty} \sum_{m\geq0} (z;q)_{N-m} \qbin{N}{m} z^m \cr
& ~=~ {1\over (z;q)_\infty} \sum_{m,n\geq0} (-1)^n q^{ {1\over2} n(n-1)}
  \qbin{N}{m} \qbin{N-m}{n} z^{m+n} \cr
&  ~=~ {1\over (z;q)_\infty}\sum_{m,n\geq0} (-1)^n q^{ {1\over2} n(n-1)}
  { (q)_N \over (q)_{m-n} (q)_n (q)_{N-m} } z^m \cr
& ~=~ {1\over (z;q)_\infty} \sum_{m\geq0} (1;q)_m { (q)_N \over (q)_{N-m}
  (q)_m } z^m 
  ~=~ {1\over (z;q)_\infty} ~=~ \sum_{M\geq0} {1\over (q)_M} z^M \cr}
$$
\Box

\vfil\eject


\footatend\vfill\eject\immediate\closeout\rfile
\baselineskip=14pt{{\bf  REFERENCES}}\bigskip{\frenchspacing%
\parindent=20pt\escapechar=` \input refs.tmp\vfill\eject}\nonfrenchspacing


\vfil\eject\end